\documentclass{article}

\usepackage{arxiv}

\usepackage[utf8]{inputenc}
\usepackage[T1]{fontenc}
\usepackage{lmodern}
\pdfmapfile{+lm.map}
\pdfmapfile{+symbols.map}
\pdfmapfile{+euler.map}

\usepackage{xcolor}
\usepackage{url}
\usepackage{booktabs}
\usepackage{amsmath,amssymb,amsthm,mathtools}
\usepackage{nicefrac}
\usepackage{microtype}
\usepackage{graphicx}
\usepackage[numbers,sort&compress]{natbib}
\usepackage{doi}
\usepackage{algorithm}
\usepackage[noend]{algpseudocode}
\usepackage{enumitem}
\usepackage{float}
\usepackage{placeins}
\usepackage{needspace}
\usepackage{array}
\usepackage{hyperref}
\usepackage{authblk}
\usepackage[capitalise,nameinlink]{cleveref}
\Crefname{appendix}{Appendix}{Appendices}
\crefname{appendix}{Appendix}{Appendices}

\makeatletter
\renewcommand\fs@ruled{%
  \def\@fs@cfont{\bfseries}\let\@fs@capt\floatc@ruled
  \def\@fs@pre{\hrule height\heavyrulewidth depth0pt\kern3pt}%
  \def\@fs@mid{\kern3pt\hrule height\lightrulewidth depth0pt\kern4pt}%
  \def\@fs@post{\kern3pt\hrule height\heavyrulewidth depth0pt\relax}%
  \let\@fs@iftopcapt\iftrue}
\makeatother
\algrenewcommand\alglinenumber[1]{\footnotesize\color{black!60}#1:}

\newtheorem{theorem}{Theorem}
\newtheorem{lemma}{Lemma}
\newtheorem{proposition}{Proposition}
\newtheorem{corollary}{Corollary}
\theoremstyle{definition}
\newtheorem{assumption}{Assumption}
\theoremstyle{remark}
\newtheorem{remark}{Remark}
\newtheorem{example}{Example}

\numberwithin{equation}{section}
\crefname{assumption}{Assumption}{Assumptions}
\crefname{algorithm}{Algorithm}{Algorithms}

\newcommand{\R}{\mathbb{R}}

\newcommand{\norm}[1]{\lVert #1\rVert}
\newcommand{\ip}[2]{\langle #1,#2\rangle}
\newcommand{\prox}{\operatorname{prox}}
\newcommand{\dom}{\operatorname{dom}}

\setlist{itemsep=2pt,topsep=4pt}

\newcommand{\modelerror}[3]{\bigl[B_f(#1,#2)-\frac{\norm{#2-#1}^2}{2#3}\bigr]_+}

\newcommand{\Cert}{\mathcal C}
\DeclareFontShape{T1}{ptm}{m}{scit}{<->ssub * ptm/m/sc}{}

\title{Universal Restart-Free Proximal Acceleration\\with Upper and Lower Curvature Adaptation}

\author[1,2]{Ilya Kuruzov\thanks{Corresponding author: \href{mailto:kuruzov.ia@phystech.edu}{kuruzov.ia@phystech.edu}}}
\author[1,2]{Alexander Vyguzov}
\author[3]{Ivan Dorofeev}
\author[1,2]{Alexander Gasnikov}
\affil[1]{Moscow Independent Research Institute of Artificial Intelligence}
\affil[2]{Innopolis University}
\affil[3]{Independent researcher}

\renewcommand{\shorttitle}{ELENCHUS: restart-free proximal acceleration}
\date{}

\hypersetup{hypertexnames=false,colorlinks=true,linkcolor=blue!50!black,citecolor=blue!50!black,urlcolor=blue!50!black,pdftitle={Universal Restart-Free Proximal Acceleration with Upper and Lower Curvature Adaptation},pdfauthor={Ilya Kuruzov, Alexander Vyguzov, Ivan Dorofeev, Alexander Gasnikov},
  pdfkeywords={universal methods, accelerated proximal gradient methods, H\"older continuity, strong convexity, parameter-free methods, restart-free acceleration}}

\usepackage{subcaption}
\begin{document}
\DeclareFontShape{T1}{lmr}{bx}{sc}{<->ssub * lmr/m/sc}{}
\maketitle

\begin{abstract}
\textsc{ELENCHUS} is a restart-free accelerated proximal-gradient method for composite convex optimization that adapts its stepsize and momentum to upper and lower curvature along a single trajectory. It requires no curvature constants, smoothness or growth exponents, or target accuracy. Backtracking uses a running energy forecast to set the upper-model allowance, while the lower controller raises its curvature estimate at a clock-capped rate and uses an affine-minorant certificate and a protected secant to lower it when refuted, preserving the momentum point and weight across each correction. When the differentiable term is smooth and strongly convex, the objective gap satisfies an accelerated linear bound at every iterate after an additive transient, with square-root dependence on the condition number. Under global smoothness, strong convexity on the trajectory ball suffices even when the global modulus vanishes. More generally, upper and lower curvature moduli on bounded regions govern the oracle cost of first reaching a target gap, recovering the optimal power-law dependence on both accuracy and the uniform-convexity constant up to logarithms under H\"older gradients and $p$-uniform convexity. Under global smoothness and without requiring global strong convexity, once all query, primal, and momentum points remain in a more strongly curved neighborhood of the minimizer, the estimate permanently recovers at least half the local modulus in finitely many steps. Experiments on composite, H\"older, quartic, and Huber objectives examine stepsize growth, two-chord protection, and local curvature recovery.
\end{abstract}


\section{Introduction}
\label{sec:intro}

We consider the composite convex optimization problem
\begin{equation}\label{problem}
  \min_{x\in\R^d}\ F(x)=f(x)+\psi(x),
\end{equation}
where $f$ is convex and continuously differentiable and $\psi$ is proper, closed and convex with a computable proximal operator. For an $L$-Lipschitz gradient and strong convexity modulus $\mu_f>0$, the upper curvature sets the stepsize and the lower curvature sets the momentum in Nesterov's $O(\sqrt{L/\mu_f}\log(1/\varepsilon))$ bound \citep{nesterov1983,nesterov2013composite,beck2009fista,nesterov2018}. Global constants can miss the geometry along a run. Unregularized logistic loss with a full-column-rank design has zero global strong-convexity modulus and positive curvature on every bounded ball. In consistent quartic regression with full column rank, the Hessian vanishes at the solution and grows without bound at infinity, yet the loss is smooth on bounded balls and has a fourth-order lower bound.

On a trial from $x$ to a proximal point $p\ne x$, the available values and gradient give the normalized tangent error
\[
  \frac{2\bigl(f(p)-f(x)-\langle\nabla f(x),p-x\rangle\bigr)}{\|p-x\|^2}.
\]
The upper model needs to hold only on the trial chord, where an additive allowance accommodates H\"older gradients even when the quotient grows on short steps. An exact lower modulus must hold in every direction, so one small secant rules out a larger modulus even though repeated large secants can miss a flat direction \citep{kuruzov2026forgetting}. Allowing slack $s_0$ changes what constitutes a refutation. A chord with secant $c<m$ remains consistent with lower curvature $m$ whenever $(m-c)\|p-x\|^2\le2s_0$. Near a quartic minimum, a short chord can meet this condition even for very small $c$.

Existing approaches address different parts of this difficulty. Universal upper models control the error of a trial, level methods obtain additional geometric information through projections, and restart schemes can exploit lower curvature by choosing the length of an accelerated stage. A single quadratic mode illustrates how the stage length affects the momentum accumulated between resets.

Consider the slow eigenmode of a quadratic with condition number $\kappa=L/\mu_f>1$. Fix the ideal parameters $\alpha=1/L$ and $\sigma=\gamma=\mu_f$ to isolate the effect of resetting $v\gets y$. With $\omega=1/\sqrt\kappa$, one step maps $(y,v)$ to $(1-\omega)(y+\omega v,v)$. Starting from $v_0=y_0$, a block of $q$ steps therefore gives
\[
 y_q=(1-\omega)^q(1+q\omega)y_0.
\]
For $q=\lceil 2\sqrt\kappa\rceil$, the multiplier is at most $(1+q\omega)e^{-q\omega}\le 3e^{-2}<1/2$. Repeating blocks of this length gives the accelerated $O(\sqrt\kappa\log(1/\varepsilon))$ iteration count on this mode.

For $q=o(\sqrt\kappa)$, the linear terms in the two factors cancel. The mean decrease in the log objective per step is $-2\log(1-\omega)-(2/q)\log(1+q\omega)\sim(q+1)/\kappa$. Thus periodic resets require asymptotically $\kappa\log(1/\varepsilon)/(q+1)$ steps to reach relative gap $\varepsilon$, taking $\varepsilon\downarrow0$ before $\kappa\to\infty$. At $q=1$ the iterate multiplier is exactly $1-1/\kappa$, the gradient descent factor. This calculation motivates separating curvature corrections from momentum resets when the estimate may change several times within one stage.

ELENCHUS (\emph{Evidence-based Lower-curvature Estimation for Nesterov's method: a Certified, Horizon-free, Universal Scheme}) revises its working estimate $\sigma_k$ along one accelerated trajectory. An observable certificate compares the aggregate of affine minorants with a running forecast, and a protected secant tests the estimate on the accepted chord. If either test fails, the controller lowers the estimate and recalibrates the minorant and forecast while carrying the accepted primal point, momentum point and separate momentum weight $\gamma_k$ into the next step. Clock-capped raises let the estimate recover when the trajectory stays in a more strongly curved region. Clipping the trial estimate at the inverse stepsize lets the first proposal grow independently of the inherited momentum weight. By default, a cached anchor supplies a second chord for the cutoff. Neither the curvature parameters nor the target accuracy enter these updates.

\subsection{Related work}

\looseness=-1 Universal methods use an additive allowance tied to the target accuracy for unknown H\"older or more general upper moduli \citep{nesterov2015universal,lan2015,nesterov2025gcb}. AGDA starts from a distance guess and adapts it for anytime H\"older guarantees, while AdaNAG and AdaAGM adapt to local smoothness \citep{ren2025agda,suh2025adaptive,wang2025adaptive}. Under strong convexity, the linear factor proved for AdaAGM depends on $\mu_f/L$ \citep{wang2025adaptive}. In ELENCHUS the forecast sets the upper allowance and also governs corrections of the lower estimate.

Sequential subspace optimization exploits curvature in selected directions, and conjugate gradients exploit the fixed operator of a positive definite quadratic \citep{narkiss2005,elad2007sesop,kuruzov2026forgetting,hestenes1952}. Catalyst and bundle-level methods instead solve auxiliary subproblems or projections \citep{lin2018catalyst,jiang2026abl}. For A-BLW, these projections yield accelerated bounds under quadratic and H\"older growth. ELENCHUS uses the original proximal operator and stores the aggregate minorant as one vector and one scalar.

\looseness=-1 Restart lengths can use strong convexity, quadratic growth or H\"olderian error bounds \citep{nesterov2013composite,necoara2019,roulet2020,liu2017}, while adaptive schedules estimate curvature or test progress \citep{lin2014,lin2015coap,fercoq2019,fercoq2020,aujol2024freefista,alamo2019,sujanani2025}. For example, SCAR tests a modulus across accelerated stages for a gradient-norm target \citep{lan2024scar}. The tuned universal restart schedule of \citet{roulet2020} uses both curvature exponents, their constants and a bound on the initial gap. Observable triggers can remain inactive even under smooth strong convexity \citep{odonoghue2015,giselsson2014,kim2018restart,su2016,moursi2023,bao2026restart,bao2026nonoccurrence}. ELENCHUS carries its momentum point and weight through a downward correction.

\looseness=-1 NAG-free corrects its trial modulus using a gradient-difference quotient and a multiplicative reduction. The resulting global rate has gradient-descent dependence on conditioning, while local acceleration needs a certified upper estimate and additional Hessian and trajectory assumptions \citep[Theorems~3 and~5]{cavalcanti2025}. Linear rates without a supplied modulus are available for accelerated forward--backward schemes \citep{li2024linear}. Li and Wang prove linear convergence of proximal acceleration for every fixed momentum parameter in $[0,1)$ \citep{liwang2026}. A$^2$GD updates both curvature estimates using a radius bound, a gradient-norm target and a decreasing positive floor \citep{xu2026adaptive}. ELENCHUS allows upward recovery and proves a global accelerated bound with an additive transient in the smooth strongly convex case.

\subsection{Contributions and organization}

ELENCHUS achieves global acceleration under unknown smoothness and strong convexity while retaining momentum through curvature corrections.

\begin{itemize}
\item \emph{Regional curvature.} Upper and lower moduli on bounded regions determine the cost of first reaching a target gap (\cref{bi:horizon}). For H\"older gradients and power-type uniform convexity, the bound captures the dependence on accuracy and $\mu_p$ without knowing either exponent or curvature constant. Smooth quartic regression gives $\widetilde O(\varepsilon^{-1/4})$ after initialization (\cref{bi:power}).
\item \emph{Restart-free acceleration.} A global upper model gives a first-attainment bound (\cref{thm:master}) with convex and mixed-curvature specializations (\cref{cor:convex,prop:mixed}). Under smooth strong convexity, every subsequent iterate satisfies an accelerated linear bound after a transient of order $\sqrt\kappa\log^2\kappa$, without supplied $L$ or $\mu_f$ (\cref{prop:smooth-energy}). Under global smoothness, strong convexity on the trajectory ball suffices.
\item \emph{Local recovery.} When the trajectory enters and stays in a region of stronger curvature, the estimate recovers a fixed fraction of the local modulus in finite time (\cref{thm:simple-local-recovery}). Local recovery also holds when the global strong-convexity modulus vanishes. The two-chord test adds at most eight accepted steps to refresh its anchor after entry.
\end{itemize}

\Cref{tab:method-matrix} compares the curvature models, inputs, and convergence guarantees of the methods discussed above.

\begin{table}[!htbp]
\centering\footnotesize
\setlength{\tabcolsep}{1.8pt}
\caption{Upper models, composite support, input requirements, rates and iteration mechanisms.}\label{tab:method-matrix}
\begin{tabular}{@{}l*{10}{c}@{}}
\toprule
& \multicolumn{2}{c}{Upper model} & & & & \multicolumn{3}{c}{Optimal order} & \multicolumn{2}{c}{Free of}\\
\cmidrule(lr){2-3}\cmidrule(lr){7-9}\cmidrule(l){10-11}
Method & GCB & H\"older & \shortstack{Proximal\\composite} & \shortstack{Parameter-\\free} & Target-free & Convex & SC & $p>2$ & Restarts & \shortstack{Auxiliary\\solves}\\
\midrule
UFGM \citep{nesterov2015universal} & \textcolor{red!70!black}{$\times$} & \textcolor{green!45!black}{$\checkmark$} & \textcolor{green!45!black}{$\checkmark$} & \textcolor{green!45!black}{$\checkmark$} & \textcolor{red!70!black}{$\times$} & \textcolor{green!45!black}{$\checkmark$} & \textcolor{red!70!black}{$\times$} & \textcolor{red!70!black}{$\times$} & \textcolor{green!45!black}{$\checkmark$} & \textcolor{green!45!black}{$\checkmark$}\\
GCB method \citep{nesterov2025gcb} & \textcolor{green!45!black}{$\checkmark$} & \textcolor{green!45!black}{$\checkmark$} & \textcolor{green!45!black}{$\checkmark$} & \textcolor{green!45!black}{$\checkmark$} & \textcolor{red!70!black}{$\times$} & \textcolor{green!45!black}{$\checkmark$} & \textcolor{red!70!black}{$\times$} & \textcolor{red!70!black}{$\times$} & \textcolor{green!45!black}{$\checkmark$} & \textcolor{green!45!black}{$\checkmark$}\\
AGDA \citep{ren2025agda} & \textcolor{red!70!black}{$\times$} & \textcolor{green!45!black}{$\checkmark$} & \textcolor{green!45!black}{$\checkmark$} & \textcolor{green!45!black}{$\checkmark$} & \textcolor{green!45!black}{$\checkmark$} & \textcolor{green!45!black}{$\checkmark$} & \textcolor{red!70!black}{$\times$} & \textcolor{red!70!black}{$\times$} & \textcolor{green!45!black}{$\checkmark$} & \textcolor{green!45!black}{$\checkmark$}\\
A-BLW \citep{jiang2026abl} & \textcolor{red!70!black}{$\times$} & \textcolor{green!45!black}{$\checkmark$} & \textcolor{red!70!black}{$\times$} & \textcolor{green!45!black}{$\checkmark$} & \textcolor{green!45!black}{$\checkmark$} & \textcolor{green!45!black}{$\checkmark$} & \textcolor{green!45!black}{$\checkmark$} & \textcolor{green!45!black}{$\checkmark$} & \textcolor{green!45!black}{$\checkmark$} & \textcolor{red!70!black}{$\times$}\\
NAG-free \citep{cavalcanti2025} & \textcolor{red!70!black}{$\times$} & \textcolor{red!70!black}{$\times$} & \textcolor{red!70!black}{$\times$} & \textcolor{green!45!black}{$\checkmark$} & \textcolor{green!45!black}{$\checkmark$} & \textcolor{red!70!black}{$\times$} & \textcolor{red!70!black}{$\times$} & \textcolor{red!70!black}{$\times$} & \textcolor{green!45!black}{$\checkmark$} & \textcolor{green!45!black}{$\checkmark$}\\
A$^2$GD \citep{xu2026adaptive} & \textcolor{red!70!black}{$\times$} & \textcolor{red!70!black}{$\times$} & \textcolor{green!45!black}{$\checkmark$} & \textcolor{green!45!black}{$\checkmark$} & \textcolor{red!70!black}{$\times$} & \textcolor{green!45!black}{$\checkmark$} & \textcolor{green!45!black}{$\checkmark$} & \textcolor{red!70!black}{$\times$} & \textcolor{green!45!black}{$\checkmark$} & \textcolor{green!45!black}{$\checkmark$}\\
AdaNAG / AdaAGM \citep{suh2025adaptive,wang2025adaptive} & \textcolor{red!70!black}{$\times$} & \textcolor{red!70!black}{$\times$} & \textcolor{red!70!black}{$\times$} & \textcolor{green!45!black}{$\checkmark$} & \textcolor{green!45!black}{$\checkmark$} & \textcolor{green!45!black}{$\checkmark$} & \textcolor{red!70!black}{$\times$} & \textcolor{red!70!black}{$\times$} & \textcolor{green!45!black}{$\checkmark$} & \textcolor{green!45!black}{$\checkmark$}\\
SCAR \citep{lan2024scar} & \textcolor{red!70!black}{$\times$} & \textcolor{red!70!black}{$\times$} & \textcolor{red!70!black}{$\times$} & \textcolor{green!45!black}{$\checkmark$} & \textcolor{red!70!black}{$\times$} & \textcolor{red!70!black}{$\times$} & \textcolor{green!45!black}{$\checkmark$} & \textcolor{red!70!black}{$\times$} & \textcolor{red!70!black}{$\times$} & \textcolor{red!70!black}{$\times$}\\
Tuned restarts \citep{roulet2020} & \textcolor{red!70!black}{$\times$} & \textcolor{green!45!black}{$\checkmark$} & \textcolor{green!45!black}{$\checkmark$} & \textcolor{red!70!black}{$\times$} & \textcolor{red!70!black}{$\times$} & \textcolor{red!70!black}{$\times$} & \textcolor{green!45!black}{$\checkmark$} & \textcolor{green!45!black}{$\checkmark$} & \textcolor{red!70!black}{$\times$} & \textcolor{green!45!black}{$\checkmark$}\\
\midrule
\textbf{ELENCHUS} & \textcolor{green!45!black}{$\checkmark$} & \textcolor{green!45!black}{$\checkmark$} & \textcolor{green!45!black}{$\checkmark$} & \textcolor{green!45!black}{$\checkmark$} & \textcolor{green!45!black}{$\checkmark$} & \textcolor{green!45!black}{$\checkmark$} & \textcolor{green!45!black}{$\checkmark$} & \textcolor{green!45!black}{$\checkmark$} & \textcolor{green!45!black}{$\checkmark$} & \textcolor{green!45!black}{$\checkmark$}\\
\bottomrule
\end{tabular}\par\smallskip
\raggedright
Parameter-free excludes supplied $L_\nu$, $\nu$, $\mu$ and $\mu_p$ or their growth counterparts.
Methods with crosses in both upper-model columns assume Lipschitz gradients. Target-free excludes prescribed final and stage tolerances.
Rate marks refer to guarantees up to logarithms. SC denotes strong convexity and $p>2$ denotes power-type curvature or growth.
NAG-free proves acceleration locally. SCAR and A$^2$GD use a gradient-norm target, and A$^2$GD also uses a radius bound.
\end{table}

\Cref{sec:problem,sec:design} introduce the curvature models and controller, \cref{sec:main,sec:additional} establish complexity and local recovery, and \cref{sec:experiments} compares gradient costs.

\FloatBarrier
\section{Problem setting and curvature model}
\label{sec:problem}

We consider the composite problem~\eqref{problem} under the following assumptions on the objective and the available proximal oracle.

\begin{assumption}\label{ass:basic}
The function $f\colon\R^d\to\R$ is convex and continuously differentiable. The function $\psi$ is
proper, closed and convex, and $\prox_{\alpha\psi}(z)=\arg\min_u\{\psi(u)+\frac{1}{2\alpha}\norm{u-z}^2\}$ can be computed for
every $\alpha>0$. The function $F=f+\psi$ has a minimizer $x^*$, and the starting point satisfies
$z_0\in\dom\psi$. Arithmetic and proximal operations are exact.
\end{assumption}

Starting from $z_0$, the method uses a trial stepsize $\widehat\alpha_0>0$ and search factor $\rho\in[1/2,1)$. Each trial calls the proximal operator once. For the analysis, fix a minimizer $x^*$ and write $F^*=F(x^*)$ and $R=\norm{z_0-x^*}$. We measure tangent error by the Bregman divergence
\[
  B_f(x,y)=f(y)-f(x)-\ip{\nabla f(x)}{y-x}.
\]

\subsection{Smoothness and upper curvature}
An $L$-Lipschitz gradient gives the quadratic upper bound $B_f(x,y)\le L\|y-x\|^2/2$. When this bound fails, a H\"older gradient of order $0<\nu<1$ still controls the tangent error by a power of the chord length, $\|y-x\|^{1+\nu}$. To cover both cases and more general growth laws, we use the \emph{global curvature bound} (GCB) of \citet{nesterov2025gcb},
\begin{equation}\label{gcbdef}
  \widehat{\mu}(t)=\sup_{\substack{\norm{x-y}\le t\\0<a<1}}
  \frac{af(x)+(1-a)f(y)-f(ax+(1-a)y)}{a(1-a)} .
\end{equation}
For $f(x)=L\|x\|^2/2$, the numerator equals $a(1-a)L\|x-y\|^2/2$. Dividing out the mixing weight leaves the same curvature scale at every point on the chord. The GCB is subadditive in $f$, so the upper moduli of individual loss terms can be added. The global specialization requires $\widehat\mu(t)<\infty$ for each finite $t$.

Taking endpoint limits in~\eqref{gcbdef} yields $B_f(x,y)\le\widehat{\mu}(\norm{x-y})$, and the condition $\widehat\mu(t)\le s$ covers chords of length at most $t$ by the allowance alone. On longer chords, scaling bounds the error by $(s/t^2)\|x-y\|^2$.

\begin{lemma}[An upper model with an allowance]\label{lem:allowance-model}
Let $f$ be convex and continuously differentiable with $\widehat{\mu}(t)<\infty$ for every $t\ge0$. Then
$a^2\widehat{\mu}(t)\le\widehat{\mu}(at)\le a\widehat{\mu}(t)$ for $t>0$ and $a\in(0,1]$. For every $s,t>0$ such that $\widehat{\mu}(t)\le s$ and all $x,y$,
\begin{equation}\label{eq:allowance-model}
 B_f(x,y)\le \frac{s}{t^2}\norm{x-y}^2+s.
\end{equation}
If $f$ is affine, the same inequality holds for every $s,t>0$.
\end{lemma}
\nopagebreak[4]
\begin{proof}
See Appendix~\ref{app:allowance}.
\end{proof}

In the stepsize convention used below, $1/\alpha\ge2s/t^2$ suffices for $B_f(x,y)\le\norm{x-y}^2/(2\alpha)+s$. This is the mechanism of universal gradient methods \citep{nesterov2015universal}, expressed directly through the original curvature bound.

\subsection{A lower modulus on the region visited by the method}
Strong convexity gives momentum a curvature scale: the graph of $f$ lies above each tangent by at least a quadratic in the distance from the point of tangency. The global lower bound with modulus $\mu_f\ge0$ reads
\begin{equation}\label{sc}
  f(u)\ge f(x)+\ip{\nabla f(x)}{u-x}+\frac{\mu_f}2\norm{u-x}^2\qquad\text{for all }u,x .
\end{equation}
We write $\kappa=L/\mu_f$ whenever $\nabla f$ is $L$-Lipschitz and $\mu_f>0$.

Even when the global modulus is zero, a positive quadratic lower bound may hold on a ball containing the accepted trajectory and $x^*$. More generally, we allow a lower model with slack on this ball $S$. The upper bound is checked on a larger ball $T$, since trial proximal points can lie farther out. For any convex $S$ of positive diameter, first take the infimum of the normalized Jensen gap over all chords of length at least $t$:
\[
 \check\mu_S(t)=\inf_{\substack{x,y\in S,\ \|x-y\|\ge t\\0<a<1}}
 \frac{af(x)+(1-a)f(y)-f(ax+(1-a)y)}{a(1-a)} .
\]
An empty infimum is $+\infty$. Now allow a slack $s>0$ and ask how large a quadratic can fit below this profile. Each length $t>0$ imposes one upper bound on its coefficient:
\[
 \frac{m}{2}t^2-s\le\check\mu_S(t)
 \quad\Longleftrightarrow\quad
 m\le\frac{2(\check\mu_S(t)+s)}{t^2}.
\]
The coefficient must satisfy all these bounds at once. Their infimum is therefore the largest admissible coefficient, which we denote by
\[
 m_S(s)=\inf_{t>0}\frac{2(\check\mu_S(t)+s)}{t^2}.
\]
We write $\widehat\mu_S$ for \eqref{gcbdef} restricted to $S$. Like the upper modulus, $m_S(s)$ uses chords throughout the region. The next lemma shows how its dependence on slack specializes under power-type uniform convexity \citep{bi:juditsky}.

\noindent\begin{minipage}{\linewidth}
\begin{lemma}[Lower quadratic model]\label{bi:model}
For a convex set $S$ of positive diameter and $s>0$, $m_S(s)$ is finite, nonnegative, nondecreasing and concave in $s$. For $x,u\in S$,
\[
 B_f(x,u)\ge \tfrac12 m_S(s)\|u-x\|^2-s,\qquad
 m_S(s)\ge 2s/\operatorname{diam}(S)^2 .
\]
If $B_f(x,u)\ge (\mu_p/p)\|u-x\|^p$ on $S$ for $p\ge2$, then
$m_S(s)\ge 2(2^{2-p}\mu_p/p)^{2/p}s^{1-2/p}$.
At $p=2$ this gives $m_S(s)\ge\mu_2$.
\end{lemma}
\nopagebreak[4]
\begin{proof}
See Appendix~\ref{bi:proof-model}.
\end{proof}\end{minipage}\par

For a smooth, strongly convex loss, the tangent error lies between two exact quadratics. Allowing an additive error on each side gives the following extension. For a common allowance $s>0$, the two models on their shared region take the form
\[
 \tfrac12 m_S(s)\|u-x\|^2-s \le B_f(x,u) \le \tfrac1{2\alpha}\|u-x\|^2+s .
\]
The proof chooses the two allowances separately, taking the lower slack on $S$ proportional to the target gap.

\subsection{Examples of upper and lower curvature}\label{sec:examples}

\begin{example}[Smooth and H\"older losses]
For least squares, $f(x)=\|Ax-b\|^2/2$ gives $B_f(x,y)=\|A(y-x)\|^2/2$. Thus $L=\|A\|^2$, and full column rank gives $\mu_f=\sigma_{\min}(A)^2$. More generally, on a convex region $T$, the bound $\|\nabla f(x)-\nabla f(y)\|\le L_\nu\|x-y\|^\nu$ for distinct $x,y$ gives
\[
 \widehat\mu_T(t)\le \frac{2^{1-\nu}}{1+\nu} L_\nu t^{1+\nu},\qquad 0\le\nu\le1 .
\]
Averaging the tangent bounds at the chord's barycenter gives this estimate. At $\nu=1$ it gives the usual $Lt^2/2$ bound, while at $\nu=0$ the assumption bounds the oscillation of the gradient. Independently, uniform convexity of order $p\ge2$ on $S$ means $B_f(x,y)\ge(\mu_p/p)\|y-x\|^p$. By \cref{bi:model} it gives an effective lower curvature of at least order $\mu_p^{2/p}s^{1-2/p}$. At $p=2$ the lower modulus is already available at zero slack.
\end{example}

\begin{example}[Quartic regression]
Consider $f(x)=\|Ax-b\|_4^4/4$ with $A\in\mathbb R^{n\times d}$ of full column rank and $b=Ax^*$. Its Hessian vanishes at $x^*$ and is bounded on every bounded ball. The upper gradient condition therefore holds there with $\nu=1$, while the global strong convexity modulus is zero. The scalar loss $u^4/4$ has Bregman divergence at least $h^4/12$ on a displacement $h$, since the difference is $\tfrac32 h^2(u+h/3)^2$. Summing over rows gives
\[
 B_f(x,y)\ge\frac{\|A(y-x)\|_4^4}{12}
 \ge\frac{\sigma_{\min}(A)^4}{12n}\|y-x\|^4.
\]
This is the pair $(\nu,p)=(1,4)$ on bounded regions. The global upper modulus is infinite at every positive chord length, which explains the role of a regional upper bound.
\end{example}

\section{Algorithm design}
\label{sec:design}

Each iteration asks two questions. Does the upper model hold on this trial chord within the allowed error? If either line-search test fails, the method shortens the step and tries again. Once the trial is accepted, the controller asks whether the evidence still supports its lower-curvature estimate.

The estimate $\sigma_k$ and the momentum weight $\gamma_k$ serve different purposes. A failed monitor test lowers the estimate for future steps and resets the minorant and forecast. The accepted point, momentum point and updated weight carry over to the next iteration.

\subsection{Upper model, curvature guess, and momentum}

We use Nesterov's momentum update with upper-model backtracking \citep{nesterov2018,nesterov2015universal}.

A flatter region can permit a longer step while the momentum weight still reflects earlier curvature. The first trial depends only on the last accepted stepsize, and clipping the lower estimate at $M$ keeps the momentum equation feasible.

Given the input points $y,v$, weight $\gamma$ and inverse trial stepsize $M=1/\alpha$, write $\sigma=\sigma_{\rm tr}=\min\{\sigma_{k-1},M\}$ in the single-trial identities below. The incoming estimate satisfies $\sigma_{k-1}\le\gamma$. The positive root of $M\theta^2=(1-\theta)\gamma+\theta\sigma$ lies in $(0,1]$, including $\theta=1$ when $\sigma=M$. The trial points and the corresponding momentum update are
\begin{equation}\label{step}
  x=y+\frac{\alpha\gamma}{\alpha\gamma+\theta}(v-y),\qquad
  p=\prox_{\alpha\psi}(x-\alpha\nabla f(x)),\qquad
  v^+=v+\frac{\alpha\sigma(x-v)-(x-p)}{\theta}.
\end{equation}
Write $d=(x-p)/\alpha$. The accepted primal point $p$ satisfies
\begin{equation}\label{eq:v-identity}
  p=(1-\theta)y+\theta v^+ .
\end{equation}
Thus $v^+=y+(p-y)/\theta$ after acceptance, and the root equation updates the weight,
\begin{equation}\label{eq:gamma-update}
  \gamma\gets\frac{\theta^2}{\alpha}=(1-\theta)\gamma+\theta \sigma ,
\end{equation}
as in Nesterov's estimate sequence. With $\sigma=\gamma=\mu_f$ and $1/\alpha=L$, this gives $\theta=1/\sqrt\kappa$. Allowing $\sigma$ to vary lets the controller revise it while preserving the momentum energy across a correction.

Let $e=\modelerror{x}{p}{\alpha}$ denote the trial error. The allowance uses the incoming energy forecast $\mathcal Z$, the observed squared radius $D$, and the accumulated momentum parameters $\chi_{k-1}=\sum_{i<k}\theta_i$. All three stay fixed during the search. The first branch reserves contraction for this error, and the second bounds its accumulation across resets:
\begin{equation}\label{allow}
  \delta=\min\Bigl\{\frac{\mathcal{Z}}8,\ \frac{D\gamma}{64(1+\chi)^2}\Bigr\}.
\end{equation}
For a fixed $\rho\in[1/2,1)$, the search starts at $\alpha_{k-1}/\rho$ and halves the step after either failed test. With $c=2B_f(x,p)/\|p-x\|^2$ on a nonzero chord it requires $M\ge\min\{c,(1+\theta/(1+\chi))\sigma_{k-1}\}$ and
\begin{equation}\label{eq:test}
  f(p)\le f(x)+\ip{\nabla f(x)}{p-x}+\frac{\norm{p-x}^2}{2\alpha}+\theta\delta .
\end{equation} The additional secant test keeps the raised estimate within $M$, so $\sigma_k\le M$ after acceptance. If the clock permits a raise above $M$, acceptance requires $c\le M$, so the upper error is zero.
Every failed trial recomputes $\theta,x,p$ from the same input state $(y,v)$. The upper model inequality holds whenever $\widehat\mu(\sqrt{2\alpha\theta\delta})\le\theta\delta$. At acceptance $e\le\theta\delta$, so the first branch of~\eqref{allow} caps the error at $\theta\mathcal Z/8$. The radius estimate charges $e$ in units of the new weight $w=\theta^2/\alpha$. The root equation and $M_k\ge\rho M_{k-1}$ give a uniform lower bound on $w/\gamma$ (Appendix~\ref{app:trajectory-state}). Since the clock advances by $\theta\le1$, these increments sum against the integrable weight $(1+\chi)^{-2}$.

The factor $\rho\in[1/2,1)$ balances the speed of tracking decreasing curvature against rejected first trials. On this range $1+\log_2(1/\rho)\le2$, and an inherited first trial preserves the bound $1/(\alpha\gamma)\le Q$ whenever $Q\ge(1-\rho)^{-2}$ (see the proof in \cref{app:ratio-progress}).

The accepted pair also supplies a global affine minorant
\begin{equation}\label{cut}
  \Lambda(z)=F(p)+\frac{\ip{x-p}{z-x}+\norm{p-x}^2/2}{\alpha}-e\le F(z).
\end{equation}
For Nesterov's energy $\mathcal{E}(y,v,\gamma)=F(y)-F^*+\frac\gamma2\norm{v-x^*}^2$, strong convexity~\eqref{sc} gives
\begin{equation}\label{energy}
  \mathcal{E}(p,v^+,\theta^2/\alpha)\le(1-\theta)\mathcal{E}(y,v,\gamma)+\frac\theta2(\sigma-\mu_f)_+\norm{x-x^*}^2+e .
\end{equation}
In this energy bound, an overestimate $\sigma>\mu_f$ contributes an additive penalty, as does upper-model error, so the step can continue after a correction without restarting. For the classical iteration with fixed $\alpha=1/L$ and $\sigma=\gamma=\mu_f$, both errors vanish and $\mathcal{E}(y_k,v_k,\mu_f)\le(1-\kappa^{-1/2})^k\mathcal{E}(y_0,v_0,\mu_f)$.

\subsection{Energy, certificate and forecast}
\label{sec:oracle}

One chord suffices for the upper test. A large secant leaves curvature in unvisited directions unknown. Since the energy~\eqref{energy} depends on $x^*$ and $F^*$, we build an observable certificate from the affine cuts~\eqref{cut}, updating their aggregate by
\begin{equation}\label{eq:aggregate}
  \ell\gets(1-\theta)\ell+\theta\Lambda.
\end{equation}
The aggregate minorizes $F$. Its regularized minimum $\min_u\{\ell(u)+\sigma\|u-y\|^2/2\}=\ell(y)-\|\nabla\ell\|^2/(2\sigma)$ gives the certificate
\begin{equation}\label{certificate}
  \Cert_\sigma(y,v,\gamma,\ell)=F(y)-\ell(y)+\frac{\norm{\nabla\ell}^2}{2\sigma}+\frac\gamma2\norm{v-y}^2.
\end{equation}

\par\noindent\begin{minipage}{\linewidth}
\begin{lemma}[Energy and certificate]\label{lem:energy-cert}
Let $f$ be $\mu_f$-strongly convex with $\mu_f>0$, let $\ell$ be affine with $\ell\le F$, and let $\sigma,\gamma>0$. Then
\[
 \begin{aligned}
  \Cert_\sigma(y,v,\gamma,\ell)&\ge \frac{\mu_f}{\mu_f+\sigma}\bigl(F(y)-F^*\bigr)+\frac\gamma2\norm{v-y}^2,\\
  \mathcal{E}(y,v,\gamma)&\le 2\Bigl(1+\frac\gamma{\mu_f}\Bigr)\Bigl(1+\frac{\sigma}{\mu_f}\Bigr)\Cert_\sigma(y,v,\gamma,\ell).
 \end{aligned}
\]
\end{lemma}
\nopagebreak[4]
\begin{proof}
See Appendix~\ref{app:energy-cert}.
\end{proof}
\end{minipage}\par

The slope term equals $(\gamma/\sigma)\|\nabla\ell\|^2/(2\gamma)$, which explains the ratio in the monitor threshold. Before attainment of the minimum, for a tentative guess $\sigma\le\mu_f$, the energy and model deficit together imply $\Cert_\sigma\le16(1+\gamma/\sigma)\mathcal Z$, so a violation refutes the guess using observable quantities. The forecast tracks the certificate scale of the current aggregate, which may stay positive even at $y=v=x^*$. Between corrections it contracts and admits the observed error:
\begin{equation}\label{eq:forecast}
  \mathcal{Z}\gets(1-7\theta/8)\mathcal{Z}+e .
\end{equation}
Since acceptance gives $e\le\theta\mathcal Z/8$, the forecast before the raise is at most $(1-3\theta/4)\mathcal Z$.

\subsection{Raising and testing the lower estimate}\label{sec:allowance}

The secant $c=2B_f(x,p)/\|p-x\|^2$ is available at acceptance. Raising toward $c$ allows recovery in stronger curvature. To limit the effect of a single stiff direction, we cap the raise by the incoming clock $\chi$. With $\gamma=\theta^2/\alpha$ and forecast~\eqref{eq:forecast}, set
\begin{equation}\label{eq:raise}
  r=\min\Bigl\{1+\frac{\theta}{1+\chi},\ \max\{1,c/\sigma\}\Bigr\},\qquad (\sigma,\gamma,\mathcal{Z})\gets r\,(\sigma,\gamma,\mathcal{Z}).
\end{equation}
The common raise preserves $\gamma/\sigma$ and increases energy by at most $r$. Its clock cap telescopes, so the product of raises is at most $1+\chi_n\le1+n$.

A small secant refutes an exact quadratic modulus above $c$ \citep{cavalcanti2025}, yet a quartic can have arbitrarily small secants and still admit a useful lower model with slack. Indeed, $(m-c)\|p-x\|^2\le2s_0$ under a lower model of curvature $m$ and slack $s_0$. A secant below $m/8$ therefore comes from a short chord. Before attainment, at a tentative $\sigma\ge m/16$, the regional certificate bound gives $\sigma\Cert_\sigma\ge m(F(p)-F^*)/18$. Dividing its square root by the chord length protects the estimate on precisely these chords. The proximal chord contributes
\begin{equation}\label{eq:fast-cutoff}
 b=\max\left\{c,\frac{\sqrt{\sigma\Cert_\sigma(p,v^+,\gamma,\ell)}}{4\|p-x\|\min\{64,1+\chi+\theta\}^{4}}\right\}.
\end{equation}
With $h=\min\{64,1+\chi+\theta\}$, rejection by the cutoff requires both $c<\sigma$ and $\Cert_\sigma<16\sigma\|p-x\|^2h^8$. Capping $h$ keeps the lower slack a positive fraction of the target. On nonstationary steps $F(p)-\ell(p)\ge\theta\|p-x\|^2/(2\alpha)>0$, so $b>0$.

The tentative estimate is corrected if the certificate exceeds its forecast bound or the estimate exceeds the protected cutoff,
\begin{equation}\label{eq:trigger}
 \Cert_\sigma>16(1+\gamma/\sigma)\mathcal Z\qquad\text{or}\qquad \sigma>b .
\end{equation}
The correction sets
\begin{equation}\label{eq:correction}
 \sigma\gets\min\{\sigma/2,b\}.
\end{equation}
The correction replaces $\ell$ by $\Lambda$ and resets the forecast. Although the old cuts remain valid, their weights reflect the discarded estimates. Rebuilding the monitor leaves both momentum variables available for the next step. Before attainment of the minimum, every tentative $\sigma\le\mu_f$ passes both tests (Appendix~\ref{app:basic}).

\subsection{Secant initialization}

The initial search requires an exact upper model. Choose any first trial $\widehat\alpha_0>0$. For a scale estimate independent of the probe stepsize, take $q=\prox_{\bar\alpha\psi}(z_0-\bar\alpha\nabla f(z_0))$ with any $\bar\alpha>0$. Return $z_0$ if $q=z_0$. Otherwise the recommended choice is
\begin{equation}\label{eq:initial-choice}
 \widehat\alpha_0=\frac{\|q-z_0\|^2}{\langle\nabla f(q)-\nabla f(z_0),q-z_0\rangle},
 \qquad \widehat\alpha_0=\bar\alpha\quad\text{if the denominator vanishes}.
\end{equation}
The accepted initial pair must satisfy
\begin{equation}\label{initial-upper}
 p_0=\prox_{\alpha_0\psi}(z_0-\alpha_0\nabla f(z_0)),\qquad f(p_0)\le f(z_0)+\ip{\nabla f(z_0)}{p_0-z_0}+\frac{\norm{p_0-z_0}^2}{2\alpha_0}.
\end{equation}
Halve $\widehat\alpha_0$ until~\eqref{initial-upper} holds, reusing $f(z_0)$ and $\nabla f(z_0)$.

Write $\mathcal M_k=(\sigma_k,\mathcal Z_k,\ell_k)$ for the committed monitor state. At $p_0=z_0$ return $z_0$ by proximal optimality. Otherwise initialize with the accepted cut $\Lambda_0$ from~\eqref{cut} at $(z_0,p_0,\alpha_0,0)$:
\begin{equation}\label{eq:initial-state}
\begin{aligned}
 y_0=v_0=y_{\mathrm{out},0}&=p_0, &\qquad D_0&=\|p_0-z_0\|^2,\\
 \chi_0&=0, & \sigma_0=\gamma_0&=1/(2\alpha_0),\\
 \mathcal M_0&=(\sigma_0,6D_0/\alpha_0,\Lambda_0).
\end{aligned}
\end{equation}

At $v_0=y_0$, the cut deficit and slope terms sum to $3D_0/(2\alpha_0)$. Since $\gamma_0=\sigma_0$, the reset multiplies this sum by four, giving $\mathcal Z_0=6D_0/\alpha_0$.
\begin{proposition}[Secant initialization]
\label{prop:initial-scale}
Under \cref{ass:basic}, the initial upper search terminates, and a stationary return is a minimizer. If $f$ has an $L$-Lipschitz gradient and is $\mu_f$-strongly convex with $\mu_f>0$, the choice~\eqref{eq:initial-choice} gives, at every nonstationary initialization,
\[
 \mu_f\le 1/\alpha_0<2L,\qquad N_{\rm init}\le2+\lceil\log_2\kappa\rceil,
\]
where $N_{\rm init}$ bounds the number of calls to each oracle in the probe and initial search. These bounds do not depend on $\bar\alpha$. Moreover, every exact initial search gives
\[
 6D_0/\alpha_0\le12\bigl(F(z_0)-F(p_0)\bigr),\qquad D_0=\norm{p_0-z_0}^2.
\]
\end{proposition}
\nopagebreak[4]
\begin{proof}
See Appendix~\ref{app:initialization}.
\end{proof}

We retain the realized initialization cost and scale in the general bounds. For a supplied passing pair, estimates that use $1/\alpha_0<2L$ require this additional scale bound.

\raggedbottom
\subsection{Two-chord cutoff and the full algorithm}\label{sec:assembling}
A short proximal step can give little information about curvature across several iterates. A cached anchor $(u,f(u))$ supplies a second displacement from the same query point. The anchor starts at $(p_0,f(p_0))$ and refreshes after every eighth accepted step. The controller takes the smaller of the two protected cutoffs for both the test~\eqref{eq:trigger} and the correction~\eqref{eq:correction}. The extra chord reuses $f(u)$ and the gradient at $x_k$. Only the proximal chord enters the raise.

\begin{figure}[H]\centering
\includegraphics[width=.78\linewidth]{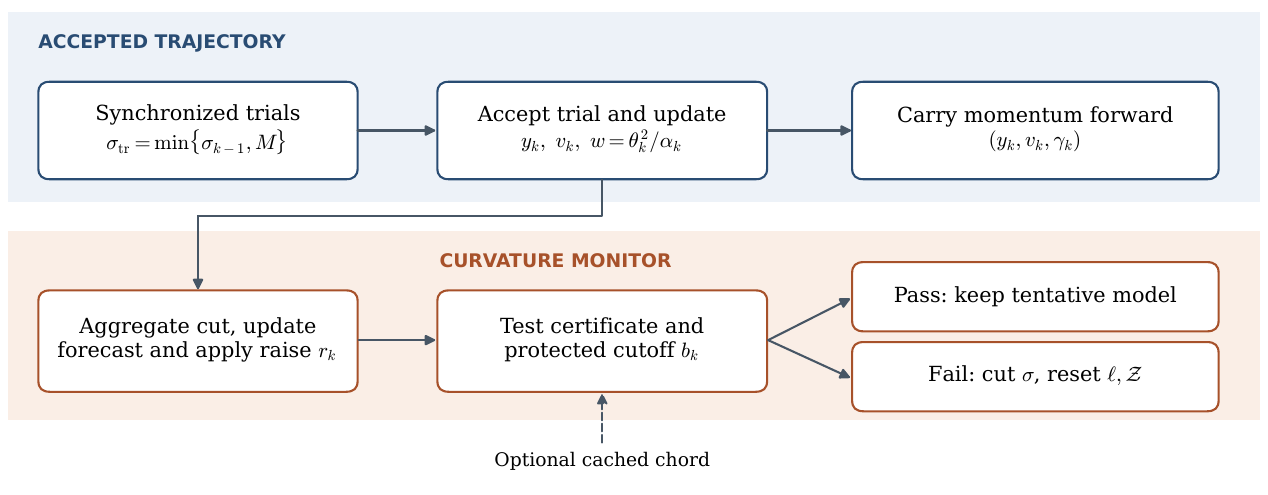}
\caption{One synchronized iteration. Failure of either trial test halves the stepsize and repeats the trial from the incoming state. After acceptance, the monitor applies the common raise $r_k$ before testing the tentative estimate. Both outcomes retain $(y_k,v_k,\gamma_k)$ with $\gamma_k=r_kw$ and pass their committed model $(\sigma_k,\mathcal Z_k,\ell_k)$ to the next iteration. The optional cached chord affects only the cutoff.}\label{fig:controller}
\end{figure}

\Cref{fig:controller} separates the accepted trajectory from the curvature monitor. The controller uses the fixed constants $8$, $16$, $64$ and $7/8$, and the output tracker $y_{\mathrm{out},k}$ compares cached values without feeding back into the iteration.

\Cref{alg:main} gives the synchronized search and calls the lower-curvature update in \cref{alg:mu-adapt} after acceptance.

\newcommand{\stepmath}[1]{\par\nobreak\vspace{1.5pt}\noindent
  \hspace*{1.2em}\(\displaystyle #1\)\par\nobreak\vspace{1.5pt}}

\begin{algorithm}[H]
\caption{\textsc{ELENCHUS} with synchronized search}
\label{alg:main}
\begin{algorithmic}[1]
\Require $f,\nabla f,\psi,\prox_{\alpha\psi}$, start $z_0\in\dom\psi$, $1/2\le\rho<1$, and a stopping rule.
\Statex Supply an initial stepsize or a positive probe stepsize $\bar\alpha$. Use two chords by default.
\State \textit{Initialize.} Cache $g_0=\nabla f(z_0)$. For a supplied initial step $\alpha$, skip the probe.
\Statex Otherwise take $q=\prox_{\bar\alpha\psi}(z_0-\bar\alpha g_0)$. If $q=z_0$, return $z_0$. Set
\stepmath{\alpha=\|q-z_0\|^2/\langle\nabla f(q)-g_0,q-z_0\rangle\quad\text{(use $\bar\alpha$ if the denominator is zero).}}
\Statex Compute $p_0=\prox_{\alpha\psi}(z_0-\alpha g_0)$. Halve $\alpha$ and recompute $p_0$ until
\stepmath{f(p_0)\le f(z_0)+\langle g_0,p_0-z_0\rangle+\|p_0-z_0\|^2/(2\alpha).}
\Statex If $p_0=z_0$, return $z_0$. Set $\alpha_0=\alpha$, $y_0=v_0=y_{\rm out,0}=p_0$ and
\stepmath{D_0=\|p_0-z_0\|^2,\quad \chi_0=0,\quad \sigma_0=\gamma_0=1/(2\alpha_0),\quad \mathcal Z_0=6D_0/\alpha_0.}
\stepmath{\ell_0(z)=F(p_0)+[\langle z_0-p_0,z-z_0\rangle+D_0/2]/\alpha_0.}
\Statex For two chords, cache $(u,f(u))=(p_0,f(p_0))$. If stopping, return $y_{\rm out,0}$. Set $k=1$.
\State \textit{Set the allowance and first trial.}
\stepmath{\delta_k=\min\left\{\frac{\mathcal Z_{k-1}}8,\frac{D_{k-1}\gamma_{k-1}}{64(1+\chi_{k-1})^2}\right\},\qquad \alpha=\alpha_{k-1}/\rho.}
\State \textit{Form a trial.} Set $M=1/\alpha$, $\sigma_{\rm tr}=\min\{\sigma_{k-1},M\}$. Take the positive root of
\stepmath{M\theta^2=(1-\theta)\gamma_{k-1}+\theta\sigma_{\rm tr}.}
\stepmath{x=y_{k-1}+\frac{\alpha\gamma_{k-1}}{\alpha\gamma_{k-1}+\theta}(v_{k-1}-y_{k-1}),\qquad p=\prox_{\alpha\psi}(x-\alpha\nabla f(x)).}
\Statex If $p=x$, return $x$.
\State \textit{Test the trial.} Compute the tangent error, secant curvature and upper-model excess:
\stepmath{B=f(p)-f(x)-\langle\nabla f(x),p-x\rangle,\quad c=2B/\|p-x\|^2,\quad e=[B-M\|p-x\|^2/2]_+.}
\Statex If either inequality below fails, set $\alpha\gets\alpha/2$ and return to Step~3.
\stepmath{e\le\theta\delta_k,\qquad M\ge\min\{c,(1+\theta/(1+\chi_{k-1}))\sigma_{k-1}\}.}
\State \textit{Accept and update.} Commit $(\alpha_k,\theta_k,x_k,y_k)=(\alpha,\theta,x,p)$ and cache $F(y_k)$. Set
\stepmath{v_k=y_{k-1}+(p-y_{k-1})/\theta_k,\qquad \overline{\mathcal Z}=(1-7\theta_k/8)\mathcal Z_{k-1}+e.}
\Statex Apply Algorithm~\ref{alg:mu-adapt} to obtain $(\gamma_k,\mathcal M_k)$ using the current anchor if enabled.
\State \textit{Update the state.} Set
\stepmath{D_k=\max\{D_{k-1},\|y_k-z_0\|^2,\|v_k-z_0\|^2\},\qquad \chi_k=\chi_{k-1}+\theta_k.}
\Statex Keep the better of $y_k$ and $y_{\rm out,k-1}$ as $y_{\rm out,k}$. For two chords, refresh the anchor to $(y_k,f(y_k))$ when $8$ divides $k$.
\Statex If stopping, return $y_{\rm out,k}$. Otherwise set $k\gets k+1$ and return to Step~2.
\end{algorithmic}
\end{algorithm}

The two iterate vectors, affine model, trial workspace and best-iterate tracker require $O(d)$ storage. Each ordinary trial uses one proximal evaluation, one gradient and at most two values of $f$, with the gradient at $z_0$ reused during initialization. On acceptance, $\psi(p)$ is evaluated once and the resulting $F(p)$ serves the cut, certificate and output comparison.

\begin{algorithm}[H]
\caption{Lower-curvature update ($\mu$-\textsc{Adapt})}\label{alg:mu-adapt}
\begin{algorithmic}[1]
\Require Monitor inputs: trial estimate $\sigma_{\rm tr}$, incoming cut $\ell_{\rm in}=\ell_{k-1}$, forecast $\overline{\mathcal Z}$, clock $\chi=\chi_{k-1}$.
\Statex Accepted trial: $(x,y,v,\alpha,\theta)=(x_k,y_k,v_k,\alpha_k,\theta_k)$, cached $F(y)$, tangent error $B$, excess $e$.
\Statex Chord data: cached $f(x),\nabla f(x)$ and, for two chords, the stored pair $(u,f(u))$.
\State \textit{Aggregate the cut.}
\stepmath{\Lambda(z)=F(y)+\frac{\langle x-y,z-x\rangle+\|y-x\|^2/2}{\alpha}-e,\qquad
\tilde\ell=(1-\theta)\ell_{\rm in}+\theta\Lambda.}
\State \textit{Raise the estimate.} Use the proximal secant $c=2B/\|y-x\|^2$ and trial weight $w=\theta^2/\alpha$.
\stepmath{r=\min\{1+\theta/(1+\chi),\max\{1,c/\sigma_{\rm tr}\}\},\qquad
(\tilde\sigma,\gamma,\widehat{\mathcal Z})=r(\sigma_{\rm tr},w,\overline{\mathcal Z}).}
\State \textit{Evaluate the certificate.}
\stepmath{\widetilde{\mathcal C}=F(y)-\tilde\ell(y)+\frac{\|\nabla\tilde\ell\|^2}{2\tilde\sigma}+\frac\gamma2\|v-y\|^2.}
\State \textit{Compute the protected cutoff.} Start with $E=\{y\}$. For two chords add $u$ if $u\ne x$. Set
\stepmath{h=\min\{64,1+\chi+\theta\},\qquad b=\min_{q\in E}\max\left\{\frac{2B_f(x,q)}{\|q-x\|^2},\frac{\sqrt{\tilde\sigma\,\widetilde{\mathcal C}}}{4h^4\|q-x\|}\right\}.}
\Statex Use $B_f(x,y)=B$ and compute $B_f(x,u)=f(u)-f(x)-\langle\nabla f(x),u-x\rangle$ from the cache.
\State \textit{Accept or correct.} If $\widetilde{\mathcal C}\le16(1+\gamma/\tilde\sigma)\widehat{\mathcal Z}$ and $\tilde\sigma\le b$, return $(\gamma,(\tilde\sigma,\widehat{\mathcal Z},\tilde\ell))$.
\Statex Otherwise set $\sigma=\min\{\tilde\sigma/2,b\}$, $\ell=\Lambda$ and $\mathcal Z=2(1+\gamma/\sigma)\Cert_\sigma(y,v,\gamma,\Lambda)$.
\Statex Return $(\gamma,(\sigma,\mathcal Z,\ell))$.
\end{algorithmic}
\end{algorithm}

To track changes in $\sigma$ through state $n$, let $U_n=\sum_{k\le n}\log r_k$ measure the raises and $J_n$ count the corrections. The additional reduction beyond halving is recorded in $V_n$, the sum of $[\log(\tilde\sigma/(2b))]_+$ over corrections, with $\tilde\sigma$ the tentative estimate and $b$ the cutoff.

\noindent\begin{minipage}{\linewidth}
\begin{proposition}[Accepted steps and safe corrections]\label{prop:basic}\label{prop:safety}
Under \cref{ass:basic}, run \cref{alg:main} with either cutoff. Every
accepted trial has $\sigma_{\rm tr}\le\theta^2/\alpha\le\gamma$, $0<\theta\le1$ and~\eqref{eq:v-identity}, and its minorant~\eqref{cut}
satisfies $\Lambda\le F$ whether or not $\sigma$ is valid. If~\eqref{sc} holds, then~\eqref{energy} holds, every
correction before attainment of the minimum happens at a tentative $\sigma>\mu_f$, and every committed estimate satisfies
$\sigma\ge\mu_f/2$. Every estimate committed by an accepted step is at most $\min\{M,b,\gamma_k\}$. Every committed state,
including the initial one and every state after a correction, passes the certificate test. For every $n$, independently of upper smoothness,
\begin{equation}\label{mastercount}
  U_n\le\log(1+\chi_n)\le\log(1+n),\qquad J_n\log2+V_n+\sum_{i=1}^n\log\frac{\sigma_{i-1}}{\sigma_{{\rm tr},i}}=\log(\sigma_0/\sigma_n)+U_n .
\end{equation}
\end{proposition}
\nopagebreak[4]
\begin{proof}
See Appendix~\ref{app:basic}.
\end{proof}\end{minipage}\par

By \eqref{mastercount}, a positive floor on $\sigma_n$ turns the logarithmic raise budget into a bound on the number of corrections without an upper smoothness constant.

\FloatBarrier

\FloatBarrier
\flushbottom
\section{Main convergence analysis}
\label{sec:main}

All results in this section apply directly to \cref{alg:main,alg:mu-adapt} with either cutoff. The trajectory bounds in \cref{lem:trajectory-state} confine the curvature estimates to a bounded region. A lower bound on clock growth then converts the regional and global clock budgets into first-attainment counts. Under smooth strong convexity, we also control each subsequent iterate.

Fix a minimizer $x^*$. Let $N_{\rm init}$ bound each oracle count during initialization and put $\log_+(t)=\max\{0,\log t\}$. For all $\mu_f\ge0$, define the analytical energy
\[
  \mathcal{P}_n=F(y_n)-F^*+\frac{\gamma_n}2\norm{v_n-y_n}^2.
\]
We take $0<\varepsilon<F(p_0)-F^*$, since $p_0$ suffices at larger targets. First-attainment bounds track the uninterrupted run, and limits refer to infinite runs. Step $k$ returns $y_k=p_k$ and the updated weights. Unindexed symbols denote a single state or trial.

\noindent\begin{minipage}{\linewidth}
\begin{lemma}[Trajectory bounds]\label{lem:trajectory-state}
Under \cref{ass:basic}, the run of \cref{alg:main} has finite initial and upper searches, and every stationary return is a minimizer. At each committed state,
\[
 D_n\le\frac{320}{7}R^2,\qquad
 \max\left\{\norm{y_n-x^*}^2,\ \norm{v_n-x^*}^2+
          \frac{2(F(y_n)-F^*)}{\gamma_n}\right\}\le\frac{132}{7}R^2.
\]
Moreover, $F(y_n)-F^*\le3\gamma_n\sqrt{D_n}(R+\sqrt{D_n})$.
\end{lemma}
\nopagebreak[4]
\begin{proof}
See Appendix~\ref{app:trajectory-state}.
\end{proof}
\end{minipage}\par

To quantify the state before the first hit, choose any valid global strong convexity modulus $\mu_f$ in~\eqref{sc}. Convexity always permits $\mu_f=0$. The preceding radius estimate supplies a lower scale even in that case. We use the larger of the two,
\begin{equation}\label{eq:effcurv}
  \mu_{f,\varepsilon}=\max\Bigl\{\mu_f,\ \frac{\varepsilon}{240R^2}\Bigr\}.
\end{equation}

Both complexity bounds count accepted steps through the clock $\chi_n=\sum_{k=1}^n\theta_k$. The following observation reduces this count to controlling the accepted inverse stepsizes.
\begin{lemma}[From a curvature ratio to clock progress]\label{lem:progress-main}
In \cref{alg:main} let $Q\ge(1-\rho)^{-2}$. If the accepted inverse stepsize satisfies $1/\alpha_k\le Q\gamma_{k-1}$, then $\theta_k>1/(2\sqrt Q)$ and the next initial trial satisfies $1/\alpha_{\rm tr}\le Q\gamma_k$. Consequently, $n$ such steps spend more than $n/(2\sqrt Q)$ units of clock for $n\ge1$.
\end{lemma}
\nopagebreak[4]
\begin{proof}
See Appendix~\ref{app:ratio-progress}.
\end{proof}

A bound on the model error over the last doubling interval suffices to stop backtracking.
\par\noindent\begin{minipage}{\linewidth}
\begin{lemma}[Acceptance on the last doubling interval]\label{lem:global-accept}
In \cref{alg:main} fix $Q\ge(1-\rho)^{-2}$ and $\lambda,s>0$. Suppose $q=2\gamma/\lambda\ge1$ and $\widehat\mu(\sqrt{8s/(Q\lambda)})\le s$. Every trial with $\max\{2\gamma,Q\gamma/2\}\le1/\alpha\le Q\gamma$ and $\theta\delta\ge s/q$ passes the upper test. Hence a search starting at ratio at most $Q$ stays within $Q$ if the allowance condition holds on that interval.
\end{lemma}
\nopagebreak[4]
\begin{proof}
See Appendix~\ref{app:allowance-test}.
\end{proof}
\end{minipage}\par

\subsection{Complexity with two curvature moduli}\label{bi:section}
Under \cref{ass:basic}, consider the run of \cref{alg:main} with $0<\varepsilon<F(p_0)-F^*$. Write $S=B(x^*,5R)$, $T=B(x^*,10R)$ and $G=\sup_{x\in S}\|\nabla f(x)-\nabla f(x^*)\|$. A lower model with slack permits a small secant on a short chord. The certificate term in the cutoff protects the estimate on precisely those chords.

\begin{lemma}[Regional curvature before attainment]\label{lem:regional-main}
Fix $\tau\ge1$ and put $s_0=\varepsilon/(16\min\{64,\tau\}^8)$ and $m=m_S(s_0)$. On any prefix before first attainment with $1+\chi_n\le\tau$, every committed state satisfies
\[
 \gamma\ge\sigma\ge m/16,\qquad \delta\ge m\varepsilon/(2^{14}\gamma\tau^2),
\]
where $\delta$ is the allowance for the next search.
\end{lemma}
\nopagebreak[4]
\begin{proof}
See Appendix~\ref{app:regional-floor}.
\end{proof}

Set $\lambda=m/16$ and $q=2\gamma/\lambda$. On a trial with $2\gamma\le1/\alpha\le Q\gamma$, the same root bound as in \cref{lem:progress-main} and the regional allowance give $q\theta\delta\ge\varepsilon/(2^{10}\sqrt Q\,\tau^2)$. This determines the tolerance $s$ below. The resulting bound uses lower curvature on $S$ and checks the upper model only on the larger ball $T$ that contains the relevant trials.

\newpage
\noindent\begin{minipage}{\linewidth}
\begin{theorem}[Two curvature moduli on bounded regions]\label{bi:horizon}
In the regional setting above, for a trial ratio $Q\ge(1-\rho)^{-2}$, set
\[
 \begin{aligned}
 \tau&=C\log^2\!\left(eQ(1+R^2/(\alpha_0\varepsilon))\right),&
 s_0&=\frac{\varepsilon}{16\min\{64,\tau\}^8},\\
 \lambda&=\frac{m_S(s_0)}{16},&
 s&=\frac{\varepsilon}{2^{10}\sqrt Q\,\tau^2}.
 \end{aligned}
\]
The lower modulus supplies the scale $\lambda$. Call $Q$ admissible if $Q\lambda\ge4G/(5R)$ and the upper modulus satisfies
\begin{equation}\label{bi:upper}
 \widehat\mu_T(t)\le Qq\lambda t^2/8+s/q
 \qquad(0<t\le20R,\ q\ge1).
\end{equation}
\begin{enumerate}[label=(\roman*),leftmargin=2em]
\item \emph{Admissibility.} An admissible $Q$ exists. For a global upper bound satisfying~\eqref{bi:upper} at every $t>0$ and $q\ge1$, the condition involving $G$ may be omitted.
\item \emph{Curvature before attainment.} For any admissible $Q$, each committed state before the first hit satisfies $\gamma\ge\sigma\ge\lambda$.
\item \emph{First attainment.} An accepted point reaches gap $\varepsilon$ or less within $\lceil2\sqrt Q\,\tau\rceil$ steps and $O(N_{\rm init}+\sqrt Q\,\tau+\log(Q\tau))$ calls to each oracle.
\end{enumerate}
One may take $C=10^3(1+C_0)$ with $C_0$ from \cref{lem:global-clock-main}.
\end{theorem}
\nopagebreak[4]
\begin{proof}
See Appendix~\ref{bi:proof-horizon}.
\end{proof}
\end{minipage}\par

The threshold $Q\lambda\ge4G/(5R)$ ensures that trials in the interval used to prove acceptance have proximal points in $T$.

\begin{corollary}[Power moduli]\label{bi:power}
Suppose $\|\nabla f(x)-\nabla f(y)\|\le L_\nu\|x-y\|^\nu$ for distinct $x,y\in T$, $0\le\nu\le1$. If $B_f(x,y)\ge(\mu_p/p)\|x-y\|^p$ on $S$ for $\mu_p>0$ and $p\ge2$, the oracle count after initialization is
\[
 \widetilde O_{p,\nu,\rho}\!\left(1+
 \min\left\{\left(\frac{L_\nu R^{1+\nu}}{\varepsilon}\right)^{2/(1+3\nu)},
 \left(\frac{L_\nu}{\mu_p^{(1+\nu)/p}}\right)^{2/(1+3\nu)}
 \varepsilon^{-2(p-1-\nu)/(p(1+3\nu))}\right\}\right).
\]
Logarithms may depend on $R^2/(\alpha_0\varepsilon)$ and $L_\nu R^{1+\nu}/\varepsilon$. For the quartic regression in \cref{sec:examples}, substituting $(\nu,p)=(1,4)$ gives $\widetilde O(\varepsilon^{-1/4})$ for fixed problem parameters. For $p=2$ and $\nu=1$ the power of $1/\varepsilon$ is zero and the bound still contains logarithms.
\end{corollary}
\nopagebreak[4]
\begin{proof}
See Appendix~\ref{bi:proof-power}.
\end{proof}

For $p>1+\nu$ the accuracy exponent matches \citet[Theorem~1]{bai2025lower}, and the dependence on $\mu_p$ matches tuned restarts up to logarithms \citep{roulet2020}. The smaller branch applies to the same run at each target. For smooth quartics the cutoff contributes a factor $\min\{64,\tau\}^2$.

\subsection{Global upper curvature and a constant lower modulus}
For a global upper model the candidate trial can be tested at any distance. A lower model valid at zero slack also yields a sharper clock budget. Write $\mu_{f,\varepsilon}$ as in~\eqref{eq:effcurv}.

The global argument first bounds the two quantities entering the allowance. The forecast controls its first branch, while the observed radius controls the second.
\begin{lemma}[State bounds before first attainment]\label{lem:global-state}
Under \cref{ass:basic}, let $\mu_f\ge0$ be a valid global modulus. Before the first hit of $0<\varepsilon<F(p_0)-F^*$ in the run of \cref{alg:main},
\[
 \gamma\ge\mu_{f,\varepsilon}/2,\qquad
 \gamma\mathcal Z\ge\mu_{f,\varepsilon}\varepsilon/96,\qquad
 D\gamma^2\ge (\mu_f^2D_0+\varepsilon^2/(180R^2))/8 .
\]
\end{lemma}
\nopagebreak[4]
\begin{proof}
See Appendix~\ref{app:prehit-state}.
\end{proof}

It remains to bound the clock while the target is still out of reach.

\par\noindent\begin{minipage}{\linewidth}
\begin{lemma}[Global clock budget before attainment]\label{lem:global-clock-main}
Under \cref{ass:basic}, take a valid global modulus $\mu_f\ge0$ and $0<\varepsilon<F(p_0)-F^*$. Suppose a prefix of \cref{alg:main} stays above the target and $1/\alpha_i\le Q\gamma_{i-1}$ at each step, with $Q\ge(1-\rho)^{-2}$. For an absolute constant $C_0$,
\[
 1+\chi_n<\frac{C_0}{2}\left[\log\frac{eR^2}{\alpha_0\varepsilon}
                  +\log^2\frac{eQ}{\alpha_0\mu_{f,\varepsilon}}\right].
\]
\end{lemma}
\nopagebreak[4]
\begin{proof}
See Appendix~\ref{app:clock-budget}.
\end{proof}
\end{minipage}\par

A global upper model permits a sharper clock budget. For fixed $Q$ and $\mu_f>0$, only the first logarithm grows as $\varepsilon\downarrow0$, giving a single logarithm of inverse accuracy. The theorem takes twice this clock bound and uses the forecast and radius estimates to keep each accepted ratio below $Q$.

\begin{theorem}[Global curvature horizon]\label{thm:master}
Let \cref{ass:basic} hold, and initialize \cref{alg:main} by~\eqref{initial-upper} from any $\widehat\alpha_0>0$.
\begin{enumerate}[label=(\roman*),leftmargin=2em]
  \item \emph{Existence of an admissible curvature ratio.} For every finite global curvature bound $\widehat{\mu}$, every valid strong convexity modulus $\mu_f\ge0$, and every $0<\varepsilon<F(p_0)-F^*$, define the \emph{model-error allowance} for each candidate ratio $Q\ge(1-\rho)^{-2}$ by
  \begin{equation}\label{eq:horizon-allowance}
 \tau(Q)=C_0\left[\log\frac{eR^2}{\alpha_0\varepsilon}
             +\log^2\frac{eQ}{\alpha_0\mu_{f,\varepsilon}}\right],\qquad
 s=\frac1{1300\sqrt Q}\min\left\{\varepsilon, \frac{\mu_f^2D_0+\varepsilon^2/(180R^2)}{\mu_{f,\varepsilon}\tau(Q)^2}\right\}.
  \end{equation}
  There exists an \emph{admissible curvature ratio} $Q\ge(1-\rho)^{-2}$ satisfying
  \begin{equation}\label{masteruniv}
    \widehat{\mu}\!\left(\sqrt{\frac{8s}{Q\mu_{f,\varepsilon}}}\right)\le s.
  \end{equation}
  \item \emph{Common horizon.} For any admissible $Q$ from part (i), an accepted point has gap at most $\varepsilon$ within
  \begin{equation}\label{eq:horizon-count}
  N=\left\lceil2\sqrt Q\,\tau(Q)\right\rceil
  \end{equation}
  accepted steps. There are at most $N_{\rm init}+T$ gradient and proximal evaluations, $N_{\rm init}+2T$ evaluations of $f$, and $N_{\rm init}+N$ evaluations of $\psi$, where
  \[
  T=2N+\log_2(Q\tau(Q)).
  \]
  
\end{enumerate}
\end{theorem}
\nopagebreak[4]
\begin{proof}
See Appendix~\ref{app:global-horizon}.
\end{proof}

Doubling $Q$ from $(1-\rho)^{-2}$ until~\eqref{masteruniv} holds requires only evaluations of the curvature bound. For Lipschitz gradients with $\mu_f>0$, take $Q=\max\{(1-\rho)^{-2},4L/\mu_f\}$. Dividing the second numerator in~\eqref{eq:horizon-allowance} by $\mu_{f,\varepsilon}$ gives $4\varepsilon/3$ for $\mu_f=0$ and at least $\mu_f D_0$ otherwise.

When $\mu_f=0$, $\mu_{f,\varepsilon}=\varepsilon/(240R^2)$ and~\eqref{masteruniv} can be solved explicitly for H\"older gradients. For an arbitrary finite GCB $\widehat{\mu}$, we measure complexity through the \emph{curvature gauge}
\begin{equation}\label{gauge}
  \mathcal N_\varepsilon=\min\bigl\{k\ge1:\ k\,\widehat{\mu}\bigl(2R/k^{3/2}\bigr)\le \varepsilon\bigr\},
\end{equation}
a reference complexity scale determined by the global curvature bound. Note that $\alpha_0^{-1}R^2/\varepsilon\ge2$ whenever $\varepsilon<F(p_0)-F^*$ by the exact initial step.

\enlargethispage{2\baselineskip}
\begin{corollary}[Convex objectives]\label{cor:convex}
Under the assumptions of \cref{thm:master}, fix $0<\varepsilon<F(p_0)-F^*$.
\begin{enumerate}[label=(\roman*),leftmargin=2em]
  \item If $\widehat{\mu}(t)\le C_\nu L_\nu t^{1+\nu}$ with $0\le\nu\le1$, the method reaches gap $\varepsilon$ within
  \[
    N_{\rm init}+O_{\nu,\rho}\Bigl(\Bigl[1+\Bigl(\frac{L_\nu R^{1+\nu}}\varepsilon\Bigr)^{\frac2{1+3\nu}}\Bigr]\log^{2+\frac{4(1-\nu)}{1+3\nu}}\frac{\alpha_0^{-1}R^2+L_\nu R^{1+\nu}}\varepsilon\Bigr)
  \]
  calls to each oracle.
  \item For every finite curvature bound the gauge~\eqref{gauge} gives a bound of $N_{\rm init}+O_\rho\bigl(\mathcal N_\varepsilon\log^6(\alpha_0^{-1}R^2\mathcal N_\varepsilon/\varepsilon)\bigr)$ calls to each oracle before first reaching gap $\varepsilon$.
\end{enumerate}
\end{corollary}
\nopagebreak[4]
\begin{proof}
See Appendix~\ref{app:convex}.
\end{proof}

For $0<\nu\le1$ the power of inverse accuracy in (i) matches the lower bound of \citet{nemirovski1983} simultaneously for all targets on one run. The logarithmic power ranges from $2$ for Lipschitz gradients to $6$ for a linear upper modulus. For a sum of $m$ H\"older terms, the gauge in part (ii) is at most the largest individual gauge at accuracy $\varepsilon/m$.

\begin{remark}[Mixed upper curvature]\label{prop:mixed}
Under \cref{ass:basic}, suppose $f$ is $\mu_f$-strongly convex, $\mu_f>0$, and
\begin{equation}\label{eq:mixed-bound}
 B_f(x,y)\le L_2\|x-y\|^2/2+L_\nu\|x-y\|^{1+\nu}/(1+\nu),\qquad 0<\nu\le1.
\end{equation}
For $0<\varepsilon<F(p_0)-F^*$ the first-hit cost after initialization is
\[
 \widetilde O_{\nu,\rho}\left(\sqrt{L_2/\mu_f}+
 (L_\nu^2/(\mu_f^{1+\nu}\varepsilon^{1-\nu}))^{1/(1+3\nu)}\right)
\]
calls per oracle. The logarithms involve the displayed terms, $(F(z_0)-F^*)/\varepsilon$ and $1/(\alpha_0\mu_f)$ (Appendix~\ref{bi:proof-power}).
\end{remark}

\subsection{Every-iterate acceleration under smooth strong convexity}\label{sec:smooth-extra}

The clock gives a bound on the energy at a specified iterate. The forecast contracts between corrections, while each correction can increase it by a factor controlled by the scale ratio.
\par\noindent\begin{minipage}{\linewidth}
\begin{lemma}[Energy on a prefix with bounded scales]\label{lem:prefix-energy}
For \cref{alg:main} let $\mu_f>0$ be a valid global modulus. Suppose $1/\alpha_i\le M$ on a prefix and $\mu_f/2\le\sigma_i\le\gamma_i\le M/2$ at every committed state, including initialization. Put $K=2M/\mu_f$. Then
\[
 \mathcal P_n\le32K\mathcal Z_0\exp(-3\chi_n/4+U_n)(2048K^3)^{J_n}.
\]
\end{lemma}
\nopagebreak[4]
\begin{proof}
See Appendix~\ref{app:shared}.
\end{proof}
\end{minipage}\par

The next theorem controls the scales and the correction count in this bound under smooth strong convexity.

\begin{theorem}[Energy and convergence of the full trajectory]\label{prop:smooth-energy}
Under \cref{ass:basic}, use the probe initialization~\eqref{eq:initial-choice}. Suppose first that $\nabla f$ is $L$-Lipschitz and $f$ has global strong-convexity modulus $\mu_f>0$. Write $\kappa=L/\mu_f$.
\begin{enumerate}[label=(\roman*),leftmargin=2em]
\item \emph{Every-iterate acceleration.} For every $0<\varepsilon<F(z_0)-F^*$,
\begin{equation}\label{polyrate}
    n\ge C\sqrt\kappa\left[\log\frac{2\bigl(F(z_0)-F^*\bigr)}{\varepsilon}+\log^2\kappa\right]
    \quad\Longrightarrow\quad \mathcal{P}_n\le\varepsilon.
  \end{equation}
The constant $C$ is absolute and independent of $\rho$ and the initial trial scale. Moreover, $\limsup_n n^{-1}\log\mathcal P_n\le-3/(8\sqrt\kappa)$, and all primal, momentum and query points converge to $x^*$.
\item \emph{Stepsizes and oracle cost.} At every committed state $1/\alpha_n\le2L$ and $\gamma_n\le2L$. Each accepted iteration $n\ge1$ has $\theta_n\ge1/(2\sqrt\kappa)$. For $n$ equal to the ceiling of the threshold in~\eqref{polyrate}, first attainment costs at most $N_{\rm init}+O(n+\log(2\kappa))$ calls per oracle.
\item \emph{Curvature on the trajectory ball.} Put $\mathcal K=\{u:\|u-x^*\|^2\le132R^2/7\}$. If $f$ has strong-convexity modulus $\mu_{\mathcal K}>0$ on this ball, \cref{prop:basic} holds with this modulus. Under global $L$-smoothness, (i) and (ii) hold with $\kappa=L/\mu_{\mathcal K}$ even when the global modulus is zero.
\end{enumerate}
\end{theorem}
\nopagebreak[4]
\begin{proof}
See Appendix~\ref{app:strong}.
\end{proof}

This covers losses that are strongly convex on compact sets, such as logistic regression with full column rank, without requiring prior knowledge of $\mathcal K$.

\begin{remark}[First hitting time versus every-iterate control]
The best-so-far output retains the target after its first hit. Under smooth strong convexity, \eqref{polyrate} also controls the current iterate. The transient $\sqrt\kappa\log^2\kappa$ accounts for corrections while learning $\mu_f$. Indeed, \eqref{mastercount} and $\sigma_n\ge\mu_f/2$ give $J_n\le\log_2(1/(\alpha_0\mu_f))+\log_2(1+n)$, and each correction costs at most a polynomial factor in $\kappa$.
\end{remark}

\section{Secants and local curvature adaptation}
\label{sec:additional}

Near $x^*$, the local strong convexity modulus $\mu_\ell$ may exceed the global modulus $\mu_f$. Because corrections only decrease $\sigma$, the upward raise rule~\eqref{eq:raise} is the sole mechanism capable of adapting to $\mu_\ell>\mu_f$.

\begin{theorem}[Recovery of the local modulus]\label{thm:simple-local-recovery}
Under \cref{ass:basic}, suppose $\nabla f$ is globally $L$-Lipschitz and consider an infinite run with probe initialization~\eqref{eq:initial-choice}. Assume $f$ is $\mu_\ell$-strongly convex for some $\mu_\ell>0$ on a convex neighborhood of $x^*$ that contains all query, primal and momentum
points from an index $k_0$ onward. In part (i), set $j=k_0$ for one chord and $j=8(\lfloor k_0/8\rfloor+1)$ for two. In the latter case the refreshed anchor first enters a test at $j+1$. Put $m_j=\min\{\sigma_j,\mu_\ell\}/2$ and $H_j=\max\{1/\alpha_j,2L\}$.
\begin{enumerate}[label=(\roman*),leftmargin=2em]
  \item If $\sigma_j\ge\mu_\ell$, set $n=j$. Otherwise there exists an index $n$ by which recovery is permanent, with
  \[
    n\le j+\left\lceil\left(\frac{2\mu_\ell}{\sigma_{j}}-1\right)(1+\chi_j)\sqrt{H_j/m_j}\right\rceil+1 .
  \]
  If $f$ has global modulus $\mu_f>0$, the factor $\sqrt{H_j/m_j}$ may be replaced by $2\sqrt{L/\mu_f}$. In both cases $\sigma_k\ge\mu_\ell/2$ for all $k\ge n$. The additional number of calls to each oracle is $O(1+n-k_0+\log_+(2L\alpha_{k_0}))$.
  \item For all $n\ge k_0$, $\sigma_n\le \sigma_{k_0}(1+\chi_n)/(1+\chi_{k_0})$. In particular $\sigma_n\ge\xi \sigma_{k_0}$ with $\xi>1$ requires $\chi_n-\chi_{k_0}\ge(\xi-1)(1+\chi_{k_0})$.
  \item For any local Lipschitz constant $L_\ell$ of $\nabla f$ on a neighborhood of $x^*$, with $0<\mu_\ell\le L_\ell$,
  $\liminf_k\sigma_k\ge\mu_\ell/2$ and
  \[
    \limsup_{n\to\infty}\frac1n\log\mathcal{P}_n\le-\frac34\sqrt{\frac{\mu_\ell}{4L_\ell}} .
  \]
In addition, $\limsup_k\alpha_k^{-1}\le2L_\ell$. If a sufficiently late committed estimate is at most $L_\ell$, then $\limsup_k\gamma_k\le L_\ell$.
\end{enumerate}
\end{theorem}
\nopagebreak[4]
\begin{proof}
See Appendix~\ref{app:simple-local-recovery}.
\end{proof}

Part (i) guarantees permanent local recovery without knowing $k_0$. The dependence on $1+\chi_{k_0}$ in part (ii) quantifies the cost of the same cap that bounds the cumulative raises by a logarithm globally. A cutoff $b_r\le L_\ell$ on the local tail supplies the entry condition for the momentum-weight bound in part (iii).

\section{Numerical experiments}\label{sec:experiments}

Every ELENCHUS run uses \cref{alg:main,alg:mu-adapt}, with two chords unless a one-chord control is specified. We compare $\rho\in\{0.5,0.95\}$ using initialization~\eqref{eq:initial-choice} and the $7/8$ forecast. Counts include initialization and rejected trials. In the tables we count calls through the first accepted $y$ with $F(y)-F_{\rm ref}\le\varepsilon(F(z_0)-F_{\rm ref})$. Bold marks the lowest counts overall and among methods using neither curvature nor accuracy inputs, including ties.

APG$^\star$ receives curvature constants, while UFGM$^\star$ receives the absolute gap target for each run and uses \citet[method~(4.1)]{nesterov2015universal}. FISTA reuses gradients at fixed trial centers. Reported runs use binary64 (Python~3.12.14, NumPy~2.3.5, one OpenBLAS thread). ELENCHUS evaluates the two acceptance inequalities and both monitor tests directly in floating point. The theory concerns exact arithmetic. The supplementary smooth logistic series uses extended precision for the loss and gradients.

\subsection{Composite objectives}

\looseness=-1 The composite suite consists of Lasso, nonnegative least squares (NNLS) and $\ell_1$-regularized logistic regression in dimensions 100 to 150. We use two conditioning levels, five seeds per level and common precomputed reference objectives.

\begin{table}[!htbp]\centering\small\setlength{\tabcolsep}{4pt}
\caption{Median gradient calls to relative gap $10^{-8}$ on composite objectives. FISTA-g and FISTA-f use gradient and function-value restart tests. Each row contains five instances.}\label{tab:bench-prox}
\begin{tabular}{@{}lrrrrrrr@{}}\toprule
 & \multicolumn{2}{c}{ELENCHUS} & \multicolumn{3}{c}{Adaptive competitors} & \multicolumn{2}{c}{Supplied data}\\
\cmidrule(lr){2-3}\cmidrule(lr){4-6}\cmidrule(l){7-8}
Family & $\rho=0.5$ & $\rho=0.95$ & FISTA & FISTA-g & FISTA-f & APG$^\star$ & UFGM$^\star$\\\midrule
Lasso, $\kappa=10^{3}$ & 87 & \textbf{39} & 98 & 57 & 65 & 107 & 25\,679\\
Lasso, $\kappa=10^{5}$ & 170 & \textbf{84} & 243 & 111 & 110 & 422 & 38\,275\\
Logistic, $\kappa=10^{3}$ & 179 & \textbf{94} & 2\,206 & 535 & 536 & 262 & 2\,571\\
Logistic, $\kappa=10^{5}$ & 322 & \textbf{176} & 4\,236 & 1\,208 & 1\,112 & 2\,344 & 6\,085\\
NNLS, $\kappa=10^{3}$ & 198 & \textbf{104} & 418 & 147 & 152 & 184 & 17\,063\\
NNLS, $\kappa=10^{5}$ & 806 & \textbf{428} & 1\,789 & 633 & 634 & 1\,767 & 34\,784\\
\bottomrule\end{tabular}\end{table}

ELENCHUS at $\rho=0.95$ uses fewer gradients than every tested competitor on each of the 30 instances (\cref{tab:bench-prox}).

\looseness=-1 The numerical supplement also reports 30 smooth ridge logistic regressions against Free-FISTA, NAG-free, AdaNAG and other baselines. ELENCHUS at $\rho=0.95$ has the smallest median gradient count at relative gap $10^{-10}$ for each of the three regularization levels.

\FloatBarrier
\subsection{H\"older and quartic regression}
We compare the methods on $\ell_1$-regularized regression with H\"older-smooth loss $f(x)=\|Ax-b\|_{1+\nu}^{1+\nu}/(1+\nu)$ and penalty $0.05\|\nabla f(0)\|_\infty\|x\|_1$. Each row uses the same three consistent Gaussian systems of size $2000\times1000$. We set the numerical target using a feasible primal reference and check its accuracy with a feasible dual bound. At $\nu=1$ the loss is quadratic.

\begin{table}[!htbp]\centering\small
\caption{Median gradient calls to relative gap $10^{-5}$ on $\ell_1$-regularized regression with $2000$ observations and $1000$ variables.}\label{tab:holder}
\begin{tabular}{@{}lrrrr@{}}\toprule
$\nu$ & ELENCHUS $\rho=0.5$ & ELENCHUS $\rho=0.95$ & FISTA & UFGM$^\star$\\\midrule
$0.1$ & 326 & \textbf{189} & 202 & 813\\
$0.2$ & 315 & \textbf{121} & 306 & 1\,563\\
$0.5$ & 71 & \textbf{43} & 48 & 619\\
$1$ & 35 & 20 & \textbf{18} & 192\\
\bottomrule\end{tabular}\end{table}

ELENCHUS gains most over FISTA at $\nu=0.2$ (\cref{tab:holder}). For the quartic $f(x)=\|Ax-b\|_4^4/4$, we draw twelve $600\times500$ Gaussian designs $A=GD$ with $G_{ij}\sim N(0,1/600)$ and diagonal scales decreasing geometrically from $1$ to $1/300$. Setting $b=G\xi$ for standard Gaussian $\xi$ gives $F^*=0$. From $z_0=0$ we use probe stepsize $1/2$ and cap runs at $300{,}000$ gradient calls (accepted steps for ELENCHUS). At $10^{-8}$ the two growth settings take about the same number of accepted two-chord steps, but $\rho=0.95$ incurs fewer rejected trials (\cref{tab:quartic}). At $\rho=0.5$, the second chord more than halves the median gradient count at both targets. FISTA exhausts its budget on all twelve systems at $10^{-10}$.

\begin{table}[!ht]\centering\small
\caption{Median gradient calls on 12 unregularized quartic systems with one or two protected chords.}\label{tab:quartic}
\begin{tabular}{@{}lrrrrrr@{}}\toprule

 &\multicolumn{2}{c}{$\rho=0.5$}&\multicolumn{2}{c}{$\rho=0.95$}&&\\
 Relative gap & One & Two & One & Two & FISTA & UFGM$^\star$\\\midrule
$10^{-8}$ & 37\,346.5 & 15\,327.5 & 8\,324.5 & \textbf{8\,274.5} & 148\,840.5 & 8\,321.5\\
$10^{-10}$ & 58\,621.5 & 23\,324 & 12\,732.5 & \textbf{12\,579.5} & $>300{,}000$ & \textbf{11\,915.5}\\
\bottomrule\end{tabular}\end{table}
\FloatBarrier
Adding $\eta\|x\|^2/2$ supplies lower curvature to the same designs. We compare with UFGM$^\star$, the closest unregularized baseline, using regularized gradients and no supplied $\eta$. FISTA exhausted its unregularized budget at $10^{-10}$ and was omitted from the ridge comparison. The Newton reference and its strong-convexity lower bound give identical first-hit counts (\cref{tab:ridge}).
Both ELENCHUS settings use fewer gradients than UFGM$^\star$ on every ridge instance at both targets. At relative gap $10^{-10}$, adding $\eta=10^{-4}$ lowers the $\rho=0.95$ median from $12{,}579.5$ to $746.5$ (a factor of $16.9$), while the UFGM$^\star$ median falls by $1.9\%$, from $11{,}915.5$ to $11{,}694.5$.
\begin{table}[!ht]\centering\small\caption{Median gradient calls after adding ridge regularization to the same 12 systems. ELENCHUS uses two chords.}\label{tab:ridge}
\begin{tabular}{@{}llrrr@{}}\toprule

$\eta$ & Relative gap & ELENCHUS $\rho=0.5$ & ELENCHUS $\rho=0.95$ & UFGM$^\star$\\\midrule
$10^{-4}$ & $10^{-8}$ & 1\,119 & \textbf{598.5} & 4\,733.5\\
 & $10^{-10}$ & 1\,391.5 & \textbf{746.5} & 11\,694.5\\
$10^{-3}$ & $10^{-8}$ & 505.5 & \textbf{261} & 1\,904.5\\
 & $10^{-10}$ & 647.5 & \textbf{328.5} & 3\,913.5\\
$10^{-2}$ & $10^{-8}$ & 204 & \textbf{104} & 567.5\\
 & $10^{-10}$ & 266 & \textbf{133} & 988\\
\bottomrule\end{tabular}\end{table}

\FloatBarrier
\Needspace{7\baselineskip}
\subsection{Curvature recovery after entry}
Consider Huber regression $f(x)=\sum_{i=1}^m h(a_i^\top x-b_i)$, where $h(u)=u^2/2$ for $|u|\le1$ and $h(u)=|u|-1/2$ otherwise. The global strong-convexity modulus is zero. On $\|Ax-b\|_\infty\le1$, the loss equals $\|Ax-b\|^2/2$. In the interior its Hessian is $A^\top A$, positive definite for our full-column-rank designs. This region supplies the local curvature in \cref{thm:simple-local-recovery}.

We draw $A$ with independent $N(0,1/m)$ entries and scale its columns geometrically from $1$ to $1/2$. On the $240\times60$ systems we also scale the first 120 rows by $0.01$. We then set $b=Ax^*$ and start at $40\xi$ ($40\times20$) or $200\xi$ ($240\times60$), with $x^*,\xi$ independent standard Gaussian vectors. With these scales fixed, we test ten and twenty fresh systems, respectively, to absolute gap $F(y)\le10^{-12}$ ($F^*=0$).

\begin{figure}[!htb]\centering
\begin{subfigure}[t]{.46\linewidth}\centering
\includegraphics[width=\linewidth]{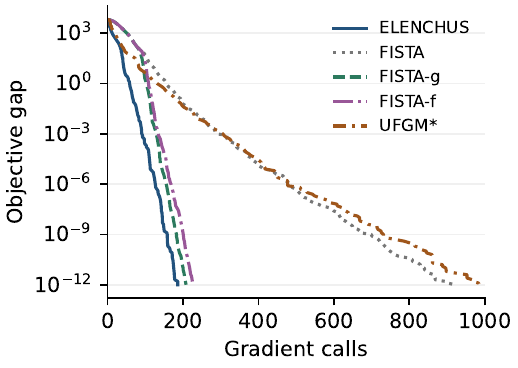}
\caption{ELENCHUS and competing methods.}\label{fig:huber-methods}
\end{subfigure}\hfill
\begin{subfigure}[t]{.46\linewidth}\centering
\includegraphics[width=\linewidth]{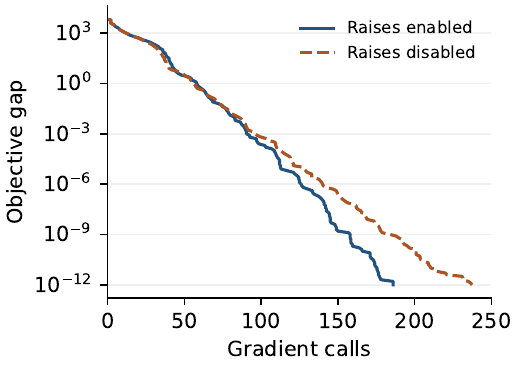}
\caption{Effect of curvature raises in ELENCHUS.}\label{fig:huber-raises}
\end{subfigure}
\caption{Huber regression on twenty $240\times60$ systems. The common starting point shows the median initial gap at zero calls. Each curve gives median gradient calls to first attain the indicated absolute gap. ELENCHUS uses $\rho=0.95$. UFGM$^\star$ receives the final absolute target $10^{-12}$.}\label{fig:huber-convergence}
\end{figure}

On the larger systems, ELENCHUS at $\rho=0.95$ has the lowest median gradient count at this target, ahead of FISTA-g and FISTA-f in Figure~\ref{fig:huber-convergence}(\subref{fig:huber-methods}). At fixed $\rho=0.95$, raises save a median of $21.3\%$ and $23.7\%$ of calls across paired runs on the two families. Figure~\ref{fig:huber-convergence}(\subref{fig:huber-raises}) shows the larger family. To isolate recovery, we replay from the first committed state after the last one with $x_k$, $y_k$, or $v_k$ outside the quadratic region. From this shared state, raises save a median of $8.4\%$ and $2.1\%$ of the remaining calls, respectively.
\FloatBarrier

\Needspace{6\baselineskip}
\section{Conclusion}
\label{sec:conclusion}

By separating the curvature estimate from the momentum weight, ELENCHUS can correct its estimate while retaining momentum. Regional curvature controls first attainment and, under global smoothness, strong convexity on the trajectory ball suffices for accelerated decay at every iterate after the transient. The synchronized search permits step doubling and enforces $\sigma_k\le1/\alpha_k$, while the cached second chord reduces the median gradient count on ill-conditioned quartics at $\rho=0.5$ without extra oracle calls per trial.

The measured gradient counts favor $\rho=0.95$ on these designs. Open questions concern the $\log^2\kappa$ factor in the strongly convex transient, the role of the finite cutoff cap in the polylogarithmic bounds, and whether recovery time can be independent of the entry clock.

\Needspace{7\baselineskip}
\section*{Declarations}

\noindent\textbf{Funding.} The research is supported by the Ministry of Science and Higher Education of the Russian Federation (Goszadaniye), project No.~FSMG-2024-0011.\par
\noindent\textbf{Competing interests.} The authors have no competing interests to declare.\par
\noindent\textbf{Data availability.} The supplement contains reference values and per-instance counts.\par
\noindent\textbf{Use of AI tools.} AI tools assisted with language editing, \LaTeX{} formatting, and figure preparation. The authors verified all mathematical derivations, code, and numerical results and take full responsibility for the content of the manuscript.

{\setlength{\bibsep}{1pt}}

\flushbottom
\appendix
\setlength{\abovedisplayskip}{6pt plus 2pt minus 2pt}
\setlength{\belowdisplayskip}{6pt plus 2pt minus 2pt}
\makeatletter
\renewcommand\subsection{\@startsection{subsection}{2}{\z@}{-1.4ex plus -.3ex minus -.2ex}{.6ex plus .1ex}{\normalfont\large\bfseries}}
\makeatother
\makeatletter
\renewcommand{\subsection}{\@startsection{subsection}{2}{\z@}%
  {-1.2ex \@plus -0.3ex \@minus -0.1ex}{0.35ex \@plus 0.15ex}%
  {\normalsize\bfseries\raggedright}}
\makeatother
\crefalias{subsection}{appendix}
\section{Proofs for Section~\ref{sec:problem}}\label{app:sec-problem}

\subsection{\texorpdfstring{Proof of \cref{lem:allowance-model}}{Proof of the allowance lemma}}\label{app:allowance}
The lower scaling inequality $a^2\widehat\mu(t)\le\widehat\mu(at)$ is \citet[Lemma 2]{nesterov2025gcb}. For the upper inequality, fix the mixing weight and expand both endpoints of a chord about their barycenter by $1/a$. The Jensen gap as a function of the dilation is convex, vanishes at zero and is unchanged in normalization. Thus the original gap is at most $a$ times the expanded gap. Taking the supremum gives $\widehat\mu(at)\le a\widehat\mu(t)$ on $\mathbb R^d$.

An endpoint limit bounds the tangent error by $\widehat\mu(\|x-y\|)$. Monotonicity handles chords no longer than $t$, and the lower scaling inequality handles the longer chords:
\[
 B_f(x,y)\le
 \begin{cases}
 s, & \|x-y\|\le t,\\
 s\|x-y\|^2/t^2, & \|x-y\|>t.
 \end{cases}
\]
For affine $f$ the divergence vanishes identically.

\subsection{Proof of Lemma~\ref{bi:model}}\label{bi:proof-model}
The endpoint limit of the Jensen quotient on $[x,u]$ is $B_f(x,u)$. Thus $\check\mu_S(\|u-x\|)\le B_f(x,u)$ on distinct pairs. The definition of $m_S$ now gives the lower model and the diameter bound. A fixed nonzero chord gives a finite upper bound for $m_S$. The infimum of affine increasing functions of $s$ is concave and nondecreasing. For the power estimate apply the two lower tangent bounds at $ax+(1-a)y$ and average them. The factor $a^{p-1}+(1-a)^{p-1}$ is at least $2^{2-p}$, so, with $c_p=2^{2-p}\mu_p/p$,

\[
 \frac{2(\check\mu_S(t)+s)}{t^2}\ge 2c_p t^{p-2}+\frac{2s}{t^2}.
\]
For $s>0$ the first term gives the claimed bound when $t\ge(s/c_p)^{1/p}$ and the second gives it when $t\le(s/c_p)^{1/p}$. At $s=0$ the first term gives the bound for $p=2$, and nonnegativity suffices for $p>2$.

\section{Proofs for Section~\ref{sec:design}}\label{app:sec-design}

The proximal optimality condition supplies three estimates used throughout the analysis. The first gives the affine cut used by the monitor, the second bounds the displacement of an accepted step, and the energy inequality shows exactly how an overestimate of lower curvature enters the momentum update. Here $\sigma$ denotes the clipped trial estimate in \cref{alg:main}.
All points and parameters refer to one accepted trial, $v$ is the input momentum point, $v^+$ the new one, $d=\alpha^{-1}(x-p)$ and $e=\modelerror{x}{p}{\alpha}$.

\begin{lemma}[Proximal inequalities]\label{basic}
Let the input primal point satisfy $y\in\dom\psi$, as every accepted state does. Every accepted trial has $\sigma\le\theta^2/\alpha\le\gamma$, $0<\theta\le1$ and $p=(1-\theta)y+\theta v^+$. Let $\mu_f\ge0$ satisfy~\eqref{sc}, where $\mu_f=0$ is always allowed. Then for all $u$
\begin{equation}\label{lower}
  F(u)\ge F(p)+\ip d{u-x}+\frac{\alpha}{2}\norm d^2-e+\frac{\mu_f}2\norm{u-x}^2 ,
\end{equation}
\begin{equation}\label{descent}
  F(p)+\frac{1}{2\alpha}\norm{p-y}^2\le F(y)+\frac\gamma2\norm{v-y}^2+e ,
\end{equation}
and $\mathcal{E}(p,v^+,\theta^2/\alpha)\le(1-\theta)\mathcal{E}(y,v,\gamma)+\frac\theta2(\sigma-\mu_f)\norm{x-x^*}^2+e$. None of this needs $\sigma\le\mu_f$.
\end{lemma}
\nopagebreak[4]
\begin{proof}
See Appendix~\ref{app:step}.
\end{proof}

\subsection{\texorpdfstring{Proof of \cref{lem:energy-cert}}{Proof of the energy lemma}}\label{app:energy-cert}
\textit{Step 1 (the certificate).} Put $\omega=\mu_f/(\mu_f+\sigma)$ and $u=(1-\omega)y+\omega x^*$. Since $\ell\le F$, completing the square and evaluating at $u$ gives
\[
 \ell(y)-\frac{\|\nabla\ell\|^2}{2\sigma}
 \le F(u)+\frac\sigma2\|u-y\|^2
 \le (1-\omega)F(y)+\omega F^* .
\]
In the last inequality the strong-convexity term $-\mu_f\omega(1-\omega)\|y-x^*\|^2/2$ cancels the quadratic penalty. Subtract from $F(y)$ and add the momentum term to obtain the first bound.

\textit{Step 2 (the energy).} The optimality of $x^*$ and strong convexity give $\norm{y-x^*}^2\le2(F(y)-F^*)/\mu_f$, and
$\norm{v-x^*}^2\le2\norm{v-y}^2+2\norm{y-x^*}^2$. Therefore
\[
  \mathcal{E}(y,v,\gamma)\le\Bigl(1+\frac{2\gamma}{\mu_f}\Bigr)\bigl(F(y)-F^*\bigr)+\gamma\norm{v-y}^2\le2\Bigl(1+\frac\gamma{\mu_f}\Bigr)\Bigl[F(y)-F^*+\frac\gamma2\norm{v-y}^2\Bigr].
\]
By Step 1 the bracket is at most $(1+\sigma/\mu_f)\Cert_\sigma(y,v,\gamma,\ell)$, since $1/\omega=1+\sigma/\mu_f\ge1$.

\subsection{\texorpdfstring{Proof of \cref{prop:basic}}{Proof of the basic proposition}}\label{app:basic}
\looseness=-1 Write $M_k=1/\alpha_k$, $w=M\theta^2=(1-\theta)\gamma+\theta \sigma_{\rm tr}$, and $a=1+\theta/(1+\chi)$. The root lies in $(0,1]$ since $\sigma_{\rm tr}\le M$. At $M<\sigma_{k-1}$ it equals one, the extra test gives $c\le M$, and the acceptance has $e=0$, $r=1$ and $v^+=p$. If the trial is unclipped, the tentative estimate is $r\sigma_{\rm tr}=\min\{a\sigma_{\rm tr},\max(\sigma_{\rm tr},c)\}\le M$ by the extra test. In both cases $\sigma_{\rm tr}\le w\le\gamma$, and a correction only decreases the estimate. The minorant and energy estimates follow from \cref{basic}. For safety we compare the certificate with auxiliary quadratics whose minimum tracks $\ell$. With $\sigma\le\mu_f$ their gap to $F$ contracts at least as fast as the forecast. The accounting identity multiplies the scalar updates of $\sigma$. We write $t=\gamma/\sigma\ge1$ and $\Cert=\Cert_\sigma(y,v,\gamma,\ell)$.

\textit{Step 1 (forecast and reset).} Averaging and resets keep $\ell\le F$, since every cut is a minorant by~\eqref{lower}. On continuation, before the raise, $\mathcal{Z}^+=(1-7\theta/8)\mathcal{Z}+e$ with $e\le\theta \mathcal{Z}/8$, so $(1-\theta)\mathcal{Z}+e\le \mathcal{Z}^+\le(1-3\theta/4)\mathcal{Z}\le e^{-3\theta/4}\mathcal{Z}$. At a reset $y=p$ and $\ell=\Lambda$, which gives
\[
  \Cert=e+\frac{\alpha}{2}\norm d^2+\frac{\norm d^2}{2\sigma}+\frac{\gamma}2\norm{v-y}^2,\qquad \mathcal{Z}=2(1+t)\Cert .
\]
The reset passes the test and gives $\mathcal Z\ge\|d\|^2/\gamma$. For $\sigma<\gamma$, the reset formula also gives
\[
 \frac{\mathcal Z}{\gamma-\sigma}\ge\frac{\|d\|^2}{\sigma^2}+\frac{2e}{\sigma}.
\]
Initialization is this reset at $y=v=p_0$ with $e=0$, $\sigma=\gamma=1/(2\alpha_0)$ and $\norm d^2=D_0/\alpha_0^2$, which gives $\mathcal{Z}_0=6D_0/\alpha_0$.

\textit{Step 2 (auxiliary quadratics).} At initialization or reset, set $W(u)=\gamma\norm{u-v-\nabla\ell/\gamma}^2/2$ and $\Phi=\ell+W$. The shift cancels the slope of $\ell$ at $v$, giving $\Phi(u)=\min\Phi+\gamma\norm{u-v}^2/2$. On continuation set $W^+=(1-\theta)W+\theta\sigma_{\rm tr}\norm{\cdot-x}^2/2$ and $\Phi^+=\ell^++W^+$. By~\eqref{step}, its Hessian is $(\theta^2/\alpha)I$ and its minimizer is $v^+$. A raise adds $(r-1)(\theta^2/\alpha)\norm{\cdot-v^+}^2/2$ to both quadratics, preserving this minimizer and $\min\Phi^+$. Thus
\begin{equation}\label{simple-deficit}
  F(y)-\min \Phi\le \mathcal{Z} .
\end{equation}
At a reset, $\min\Phi=\ell(v)+\|d\|^2/(2\gamma)$. Completing the square gives
\[
\begin{aligned}
 F(y)-\min\Phi
 &=e+\frac\alpha2\|d\|^2-\langle d,v-y\rangle-\frac{\|d\|^2}{2\gamma}\,\\
 &\le e+\frac\alpha2\|d\|^2+\frac\gamma2\|v-y\|^2\le\mathcal Z .
\end{aligned}
\]
For a continuing step, write $w=\theta^2/\alpha$ and let $q$ be the weighted mean from the proof of \cref{basic}. Dropping the nonnegative weighted variance and minimizing the remaining quadratic gives
\[
 \min \Phi^+\ge(1-\theta)\min \Phi+\theta(F(p)-e)
       +\theta\langle d,q-x\rangle-(1-\theta)\alpha\|d\|^2/2.
\]
Here the minimization contributes $-\theta^2\|d\|^2/(2w)=-\alpha\|d\|^2/2$. Since $\theta(q-x)=(1-\theta)(x-y)$, the cut at $y$ bounds the last two terms below by $(1-\theta)(F(p)-F(y)-e)$. This yields $F(p)-\min \Phi^+\le(1-\theta)(F(y)-\min \Phi)+e\le\mathcal Z^+$. A raise keeps $\min \Phi$ and increases $\mathcal{Z}$.

\textit{Step 3 (stage induction).} Write $\mathcal E=F(y)-F^*+\gamma\|v-x^*\|^2/2$. We give the induction for the forecast $\mathcal Z^+=(1-7\theta/8)\mathcal Z+e$ before the common raise. Suppose the lower model has slack $s_0$, the stage starts with $\mathcal E\le\eta\mathcal Z$, $W(x^*)\le8\eta\mathcal Z$, and $\eta\ge1$. Throughout the stage assume $\sigma\|y-x^*\|^2/2\le F(y)-F^*$ and $s_0\le\eta\mathcal Z/8$. Convexity of squared distance on $[y,v]$ then gives $\sigma\|x-x^*\|^2/2\le\mathcal E$ since $\sigma\le\gamma$. The two recursions give
\[
 \mathcal E^+\le (1-\theta)\eta\mathcal Z+e+\theta s_0\le\eta\mathcal Z^+,
 \qquad W^+(x^*)\le[8(1-\theta)+\theta]\eta\mathcal Z\le8\eta\mathcal Z^+.
\]
A raise multiplies the forecast by $r$ and adds at most $(r-1)\mathcal E^+$ to either left side. The deficit bound~\eqref{simple-deficit} uses $1-\theta\le1-7\theta/8$.

On a strict clip $M<\sigma_{k-1}$, acceptance gives $c\le M$ and $e=0$. For a lower model of modulus $m>0$ and slack $s_0$, evaluate the proximal lower inequality at $x^*$ and maximize $\langle d,x-x^*\rangle-m\|x-x^*\|^2/2$. This gives
\[
 F(p)-F^*\le s_0+\bigl(1/(2m)-1/(2M)\bigr)\|d\|^2.
\]
Thus $F(p)-F^*>s_0$ implies $M>m$. Globally use $m=\mu_f>0$, $s_0=0$ and $\eta=1$. Between corrections an unclipped estimate can only increase. Consequently a pre-attainment stage ending at a tentative $\sigma\le\mu_f$ is clip-free and all its guesses satisfy this bound. At a reset, completing the square in~\eqref{lower} gives $F(p)-F^*\le\|d\|^2/(2\sigma)+e$. Strong convexity and the reset certificate yield $\mathcal E\le(1+2t)(F(p)-F^*)+\gamma\|v-y\|^2\le\mathcal Z$ and $W(x^*)\le2\mathcal E+\|d\|^2/\gamma\le3\mathcal Z$, including at initialization.

\textit{Step 4 (certificate).} Write $W=w_0+\gamma\|\cdot-u_W\|^2/2$ with $w_0\ge0$. Since $v$ minimizes $\Phi$, the certificate equals $F(y)-\min\Phi+W(y)+t\gamma\|v-u_W\|^2/2$. Split $(y-u_W,\sqrt t(v-u_W))$ at $x^*$. The triangle inequality bounds the last two terms by $w_0+(\sqrt{t\mathcal E}+\sqrt{(1+t)(W(x^*)-w_0)})^2$. With $u=t/(1+t)$, discarding the net nonpositive contribution of $w_0$ gives
\[
 \frac{\Cert}{(1+t)\mathcal Z}\le1-u+\eta(\sqrt u+\sqrt8)^2
 \le\eta(1+\sqrt8)^2<16\qquad(1\le\eta\le16/15).
\]
The middle expression increases on $[0,1]$.

\textit{Step 5 (cutoff and lower bound).} On a nonstationary step, $F(p)-\ell(p)>0$, so $b>0$. Continuation requires $\sigma\le b$, and a correction returns a value at most $b$. Under \eqref{sc}, every tested chord has secant at least $\mu_f$, so their minimum protected cutoff is also at least $\mu_f$. The certificate bound in Step 4 rules out a correction at a pre-attainment tentative estimate at most $\mu_f$. A strict clip at $M=\mu_f$ reaches a minimizer with $p=v^+$, and even a correction then retains $\mu_f/2$. More generally, if the current step attains the minimum without such a clip, the induction in Step 3 still applies to a stage ending at a valid guess. The exact initial upper test gives $1/\alpha_0\ge\mu_f$. Induction therefore gives the claimed floor.

\textit{Step 6 (accounting).} By~\eqref{eq:raise}, $\log r_k\le\log(1+\theta_k/(1+\chi_{k-1}))=\log\frac{1+\chi_k}{1+\chi_{k-1}}$, so $U_n\le\log(1+\chi_n)\le\log(1+n)$. A correction with tentative estimate $\tilde\sigma$ returns $\tilde\sigma\exp(-\log2-[\log(\tilde\sigma/(2b))]_+)$, and raises multiply $\sigma$ by $r_k$. Each trial clip contributes the factor $\sigma_{{\rm tr},i}/\sigma_{i-1}\le1$. Multiplying these factors, raises and corrections gives~\eqref{mastercount}.

\subsection{Proof of Proposition~\ref{prop:initial-scale}}\label{app:initialization}

For any positive seed the exact search terminates by the displacement argument of Appendix~\ref{app:trajectory-state}. For the probe in~\eqref{eq:initial-choice}, $q=z_0$ implies stationarity by proximal optimality. Otherwise its secant $c$ is nonnegative since $c\|q-z_0\|^2=B_f(z_0,q)+B_f(q,z_0)$. At $c=0$ both terms vanish and the original probe passes.

Under smooth strong convexity $\mu_f\le c\le L$. The search starts at inverse scale $c$, doubles only below $L$, and every accepted scale is at least $\mu_f$ by the exact test. Hence $\mu_f\le1/\alpha_0<2L$. With at most $1+\lceil\log_2(L/c)\rceil$ trials after the probe, reuse of $f(z_0),\nabla f(z_0)$ gives the claimed count (two gradients and one final value of $\psi$). For any passing pair, proximal strong convexity and the exact upper test give $D_0/(2\alpha_0)\le F(z_0)-F(p_0)$, proving the forecast bound.

\subsection{Proof of Lemma~\ref{basic} (proximal inequalities)}\label{app:step}
All three inequalities add the prox optimality condition to the tangent and upper model inequalities of $f$, and the last one averages the first at $y$ and $x^*$.

\textit{Step 1 (scalars and points).} The root equation places $\theta^2/\alpha$ between $\sigma$ and $\gamma$. At $\theta=1$ the root polynomial has value $M-\sigma\ge0$, hence its positive root is at most one. Multiply the query equation by $\theta$ and use the root equation to obtain
\[
 x=(1-\theta)y+\theta q,
 \qquad
 q=\frac{(1-\theta)\gamma v+\theta\sigma x}{\theta^2/\alpha}.
\]
By~\eqref{step}, $v^+=q-\alpha d/\theta$. Since $p=x-\alpha d$, the displayed identity gives $p=(1-\theta)y+\theta v^+$.

\textit{Step 2 (the minorant).} Add the subgradient inequality from $d-\nabla f(x)\in\partial\psi(p)$ to the lower model at $x$. Using the upper model at $p$ gives $F(u)\ge F(p)+\ip d{u-p}-\alpha\norm d^2/2-e+\mu_f\norm{u-x}^2/2$. Now $x-p=\alpha d$ gives~\eqref{lower}.

\textit{Step 3 (descent).} Strong convexity of the proximal objective bounds its value at $p$ by its value at $y$ minus $\norm{y-p}^2/(2\alpha)$. The upper model at $p$ and the tangent bound at $y$ then give
\[
 F(p)+\frac{\norm{p-y}^2}{2\alpha}
 \le F(y)+\frac{\norm{x-y}^2}{2\alpha}+e .
\]
Put $t=M/\gamma$. From $t\theta^2\ge1-\theta$ we get $1+t\theta\ge(1+\sqrt{1+4t})/2\ge\sqrt t$. The query formula now gives $\norm{x-y}^2/\alpha\le\gamma\norm{v-y}^2$, proving~\eqref{descent}.

\textit{Step 4 (energy).} The combination of~\eqref{lower} at $y$ with weight $1-\theta$ and at $x^*$ with weight $\theta$ is
\[
  F(p)\le(1-\theta)F(y)+\theta F^*-\ip d{(1-\theta)y+\theta x^*-x}-\frac{\alpha}{2}\norm d^2+e-\frac{\theta\mu_f}2\norm{x-x^*}^2 ,
\]
after dropping $-\frac{(1-\theta)\mu_f}2\norm{y-x}^2$. Since $q$ is the weighted mean of $v$ and $x$ with weights $(1-\theta)\gamma$ and $\theta \sigma$,
\[
  \frac{\theta^2}{2\alpha}\norm{v^+-x^*}^2\le\frac{(1-\theta)\gamma}2\norm{v-x^*}^2+\frac{\theta \sigma}2\norm{x-x^*}^2-\theta\ip d{q-x^*}+\frac{\alpha}{2}\norm d^2.
\]
Adding the two displays, the terms in $\norm d^2$ cancel, and so do the linear terms, by $x=(1-\theta)y+\theta q$.

\section{Proofs for Section~\ref{sec:main}}\label{app:sec-main}

\subsection{Proof of Lemma~\ref{lem:trajectory-state}}\label{app:trajectory-state}

\textit{Step 1 (initialization and finite searches).}
The initial cut gives $2\alpha_0(F(p_0)-F^*)+\|p_0-x^*\|^2\le R^2$. Thus $D_0\le4R^2$, $\|p_0-x^*\|^2+2(F(p_0)-F^*)/\gamma_0\le2R^2$, $\mathcal Z_0\le24R^2/\alpha_0$, and $\varepsilon<R^2/(2\alpha_0)$. The initial upper test gives $\mu_f\le1/\alpha_0$ on a nonzero chord, hence $\mu_{f,\varepsilon}\le1/\alpha_0$. A zero chord certifies optimality.

For the initial trials let $r_\alpha=\|p_\alpha-z_0\|>0$. Monotonicity of $\partial\psi$ gives $(r_\alpha-r_{\alpha_*})(r_\alpha/\alpha-r_{\alpha_*}/\alpha_*)\le0$ for a fixed trial $\alpha_*$. Thus $r_\alpha/\alpha\ge r_{\alpha_*}/\alpha_*>0$ for $\alpha\le\alpha_*$. Comparison of the proximal objective with its value at $z_0$, using an affine minorant of $\psi$, gives $r_\alpha\to0$, so differentiability gives $B_f(z_0,p_\alpha)/r_\alpha\to0$ and the exact test passes. On ordinary trials $\theta\asymp\sqrt\alpha$ and $\|x-y\|=O(\sqrt\alpha)$. The proximal comparison at $y$ similarly gives $\|p-x\|=O(\sqrt\alpha)$. Continuity of $\nabla f$ gives $B_f(x,p)=o(\sqrt\alpha)\le\theta\delta$ eventually. The second test is automatic once $M\ge2\gamma$, since the prospective raise is at most $2\sigma_{k-1}$. The allowance is positive since $D\ge D_0>0$ and a continuation or a nonstationary reset keeps $\mathcal Z>0$.

\textit{Step 2 (radius and normalized errors).}
Write $w_i=M_i\theta_i^2$. The root equation gives $w_i/\gamma_{i-1}\ge1-\theta_i$. The first accepted inverse scale satisfies $M_1\ge\rho M_0\ge\gamma_0$. On subsequent steps the search and raise cap give
\[
 M_i\ge M_{i-1}/2,
 \qquad
 \gamma_{i-1}=r_{i-1}M_{i-1}\theta_{i-1}^2\le2M_{i-1}\theta_{i-1}^2.
\]
Hence $M_i\ge\gamma_{i-1}/4$ and $w_i/\gamma_{i-1}\ge\max\{1-\theta_i,\theta_i^2/4\}\ge d_0=3-2\sqrt2$. When $i=1$ or $\theta_{i-1}\le1/2$, $M_i\ge\gamma_{i-1}$ and the same minimum with $\theta_i^2$ gives $d_1=(3-\sqrt5)/2$. Otherwise $\chi_{i-1}>1/2$. By~\eqref{allow},
\[
 \frac{2e_i}{w_i}\le\frac{D_n\gamma_{i-1}}{32w_i}
 \left(1+\frac{\theta_i}{1+\chi_{i-1}}\right)
 \int_{\chi_{i-1}}^{\chi_i}\frac{dt}{(1+t)^2}.
\]
The parenthesis is at most $2$ in the $d_1$ case and $5/3$ otherwise. The integrals sum to less than one, so the normalized error sum is at most $(D_n/32)\max\{2/d_1,5/(3d_0)\}<13D_n/40$ for $D_n>0$. Write $S_n=\sum_{i\le n}2e_i/w_i$.

Divide the energy recursion by $w_i/2$. The weights $(1-\theta_i)\gamma_{i-1}/w_i$ and $\theta_i\sigma_{{\rm tr},i}/w_i$ sum to one. Thus before the raise the new normalized energy is at most the weighted average of the old one and $\|x_i-x^*\|^2$, plus $2e_i/w_i$. A raise can only decrease the normalized gap.

Since $x_i\in[y_{i-1},v_{i-1}]$ and $p_i\in[y_{i-1},v_i]$, the maximum $T_n$ of the two terms in the trajectory bound satisfies $T_n\le4R^2+S_n$. Thus \[ D_n\le\max\{D_0,2R^2+2\max_{j\le n}T_j\}\le16R^2+2S_n\le16R^2+\tfrac{13}{20}D_n,\] yielding $D_n\le320R^2/7$ and $T_n\le132R^2/7$.

For the observed-radius estimate, subtract $R^2$ from the normalized energy and call the result $\mathcal Q_i$. The initial cut gives $\mathcal Q_0\le2R\sqrt{D_0}$. In the preceding recursion, the shifted query term is $\|x_i-z_0\|^2-2\langle x_i-z_0,x^*-z_0\rangle\le D_n+2R\sqrt{D_n}$ for $i\le n$. Therefore $\mathcal Q_n\le D_n+2R\sqrt{D_n}+S_n$. Expressing $\|v_n-x^*\|^2-R^2$ in the same coordinates and using $\|v_n-z_0\|\le\sqrt{D_n}$ gives $F(y_n)-F^*\le3\gamma_n\sqrt{D_n}(R+\sqrt{D_n})$.

\subsection{Proof of Lemma~\ref{lem:progress-main}}\label{app:ratio-progress}
\textit{Step 1 (step size and clock).} Fix $Q\ge(1-\rho)^{-2}\ge4$. For $M/\gamma\le Q$, the root equation gives $Q\theta^2\ge1-\theta$, hence
\[
 \sqrt Q\,\theta\ge\frac{\sqrt{1+4Q}-1}{2\sqrt Q}\ge\frac{\sqrt{17}-1}{4}>\frac12 .
\]

\textit{Step 2 (inheritance).}
Write $Q_{\rm old}=1/(\alpha_{\rm old}\gamma_{\rm old})$ for the previous accepted ratio and $Q_*=\max\{Q_{\rm old},(1-\rho)^{-2}\}$. Its step parameter satisfies $Q_{\rm old}\theta_{\rm old}^2\ge1-\theta_{\rm old}$. If $\theta_{\rm old}<\sqrt{\rho/Q_*}\le1-\rho$, this would give $Q_{\rm old}\theta_{\rm old}^2>\rho$ and $Q_*\theta_{\rm old}^2<\rho$ simultaneously. Thus $\theta_{\rm old}^2\ge\rho/Q_*$. The new weight satisfies $\gamma\ge\theta_{\rm old}^2/\alpha_{\rm old}$, including after a correction. The first trial therefore obeys
\begin{equation}\label{simp:inherit}
 \frac{1}{\alpha_{\rm tr}\gamma}
 =\frac{\rho}{\alpha_{\rm old}\gamma}\le\frac{\rho}{\theta_{\rm old}^2}\le Q_* .
\end{equation}

With $b_i$ rejected trials, $M_i=2^{b_i}\rho M_{i-1}$ exactly. Taking logarithms and summing gives $\sum_{i=1}^n b_i=n\log_2(1/\rho)+\log_2(\alpha_0/\alpha_n)$.

\subsection{Proof of Lemma~\ref{lem:global-accept}}\label{app:allowance-test}
We first check the upper test on the last doubling interval. Fix $q=2\gamma/\lambda\ge1$ and $t=\sqrt{8s/(Q\lambda)}$ as in the statement. Since $\widehat\mu(t)\le s$, chord scaling gives $\widehat\mu(t/q)\le s/q$. The quadratic coefficient in~\eqref{eq:allowance-model} at scale $t/q$ is $sq/t^2=Q\gamma/4$. On the stated interval, $M\ge Q\gamma/2$ and $\theta\delta\ge s/q$, so
\[
 B_f(x,p)\le \tfrac14 Q\gamma\|p-x\|^2+s/q
           \le \tfrac12 M\|p-x\|^2+\theta\delta .
\]

The synchronization test also passes there because $M\ge2\gamma\ge(1+\theta/(1+\chi))\sigma_{k-1}$. Since $Q\ge4$, the interval begins at $Q\gamma/2$. A search starting below it can at most double across its lower endpoint, and a search starting inside it passes immediately. In either case the accepted scale is at most $Q\gamma$.

\subsection{Proof of Lemma~\ref{lem:regional-main}}\label{app:regional-floor}
\textit{Step 1 (region and certificate).}
By \cref{lem:trajectory-state}, the accepted points lie in $S$, as do the trial queries $x\in[y,v]$ formed from the preceding state. Write $m=m_S(s_0)$. At a point $y\in S$ with $F(y)-F^*>\varepsilon$, optimality at $x^*$ and the lower model give $m\|y-x^*\|^2/2\le F(y)-F^*+s_0$. For any affine minorant $\ell\le F$, evaluating its regularized model at $x^*$ therefore gives
\[
 \Cert_{m/16}\ge F(y)-F^*-\frac m{32}\|y-x^*\|^2
 \ge F(y)-F^*-\frac{F(y)-F^*+s_0}{16}\ge\frac89\bigl(F(y)-F^*\bigr),
\]
since $s_0\le\varepsilon/16$. For fixed $y,v,\gamma,\ell$, the product $\beta\Cert_\beta$ is affine in $\beta$ with nonnegative slope by~\eqref{certificate}. Hence
\begin{equation}\label{eq:regional-certificate-floor}
 \beta\Cert_\beta\ge\frac{m}{18}\bigl(F(y)-F^*\bigr)\qquad(\beta\ge m/16).
\end{equation}
Apply this bound at the accepted point $p$ while $F(p)-F^*>\varepsilon$, using the induction hypothesis for the incoming floor. Without clipping $\tilde\sigma\ge\sigma_{k-1}\ge m/16$. On a strict clip the strict-clip bound in Step~3 of Appendix~\ref{app:basic} gives $\tilde\sigma=M>m$.
For each tested endpoint $q$ (the proximal point or the cached earlier iterate) apply the lower model on $S$. Write $c_q=2B_f(x,q)/\|q-x\|^2$. If $c_q\ge m/8$ then its protected cutoff is at least $m/8$. Otherwise $(m-c_q)\|q-x\|^2\le2s_0$ gives $\|q-x\|^2\le\varepsilon/(7m\min\{64,\tau\}^8)$. Put $h_k=\min\{64,1+\chi_k\}$ and $h_*=\min\{64,\tau\}$. On a prefix with $1+\chi_k\le\tau$ we have $h_k\le h_*$, so
\[
 \frac{\tilde\sigma\Cert_{\tilde\sigma}}{16\|q-x\|^2 h_k^8}
 \ge \frac{m(F(p)-F^*)}{18}\frac{7m h_*^8}{16\varepsilon h_k^8}
 > \frac{7m^2}{288}>\frac{m^2}{64}.
\] Thus every candidate is at least $m/8$, as is their minimum.

\textit{Step 2 (safety on a stage).}
The strict-clip bound in Step 3 of Appendix~\ref{app:basic} gives $M>m$ before attainment. A stage ending at a tentative guess at most $m/2$ is therefore clip-free, with all earlier guesses at most $m/2$. To apply that induction with $\eta=16/15$, it remains to check the reset and slack bounds.

At the reset, evaluate the cut at $x^*$ with the lower term $m\|x-x^*\|^2/2-s_0$ and complete the square with $d=(x-p)/\alpha$. Since $s_0<(F(p)-F^*)/16$ and $\sigma\le m/2$,
\[
 F(p)-F^*\le \frac{\|d\|^2}{2m}+e+s_0
 \quad\Longrightarrow\quad
 F(p)-F^*\le\frac{16}{15}\left(\frac{\|d\|^2}{2\sigma}+e\right).
\]
The lower model at $x^*$ also yields $\sigma\|p-x^*\|^2/2\le(F(p)-F^*+s_0)/2\le F(p)-F^*$. In the reset calculation in Step 3 of Appendix~\ref{app:basic}, these two estimates bound the energy by $(16/15)\mathcal Z$ and $W(x^*)$ by $(32/15)\mathcal Z+\|d\|^2/\gamma\le(16/5)\mathcal Z$. The calculation covers initialization if its guess is at most $m/2$. A larger initial guess can enter this range only after a correction.

The pre-hit gap gives $s_0<\mathcal E/16\le(16/15)\mathcal Z/16$ at each induction step. Steps 3--4 of Appendix~\ref{app:basic} therefore preserve the energy and certificate bounds with $\eta=16/15$, since $s_0<\eta\mathcal Z/16\le\eta\mathcal Z/8$.

A certificate rejection before the hit therefore requires a tentative guess greater than $m/2$. Rejection by the cutoff requires a tentative guess greater than $b\ge m/8$ under the clock bound of Step 1. On either rejection the new estimate is at least $m/16$.

\textit{Step 3 (initialization).}
The exact initial proximal inequality gives $\varepsilon<F(p_0)-F^*\le R^2/(2\alpha_0)$. If $\|p_0-z_0\|<R/2$, then $\|p_0-x^*\|>R/2$, so
$mR^2/8 < F(p_0)-F^*+s_0 \le \tfrac{17}{16}\bigl(F(p_0)-F^*\bigr).$
Otherwise the exact initial upper test gives $m\le1/\alpha_0+8s_0/R^2<5/(4\alpha_0)$. Both cases give $\sigma_0=1/(2\alpha_0)\ge m/16$. Raises preserve this floor. A strict pre-hit clip has $M>m$ before correction, and Steps 1 and 2 propagate it on a pre-hit prefix within the clock bound.

\textit{Step 4 (paying for the tolerance).}
At a committed state before the hit, use the passing test and~\eqref{eq:regional-certificate-floor} to obtain
\[
 \frac{m(F(y)-F^*)}{18}\le\sigma\Cert_\sigma
 \le16(\sigma+\gamma)\mathcal Z\le32\gamma\mathcal Z .
\]
For the radius branch, $R\le\sqrt D+\|y-x^*\|$ and $\|y-x^*\|\le2\sqrt{(F(y)-F^*)/m}$. Apply the observed-radius bound of \cref{lem:trajectory-state} with $X=D\gamma^2/m$. Using $\gamma\ge m/16$ gives
\[
 \begin{aligned}
 F(y)-F^*&\le6\gamma D+6\sqrt{X(F(y)-F^*)}\\
 &\le96X+6\sqrt{X(F(y)-F^*)}\\
 &\le114X+\tfrac12\bigl(F(y)-F^*\bigr).
 \end{aligned}
\]
The last step uses $6\sqrt{X(F(y)-F^*)}\le18X+(F(y)-F^*)/2$. Thus $D\gamma^2\ge m(F(y)-F^*)/256$. Substitution in the two branches of~\eqref{allow} yields $\delta\ge m\varepsilon/(2^{14}\gamma\tau^2)$ whenever $1+\chi\le\tau$. For $1/\alpha\le Q\gamma$, \cref{app:ratio-progress} gives $\theta>1/(2\sqrt Q)$. With $q=2\gamma/\lambda=32\gamma/m$ this gives
$q\theta\delta\ge\frac{\varepsilon}{2^{10}\sqrt Q\,\tau^2}=s .$

\subsection{Proof of Lemma~\ref{lem:global-state}}\label{app:prehit-state}

A pre-hit prefix has $F(y_k)-F^*>\varepsilon$ at every committed point, including $p_0$. All radius estimates used below follow from \cref{lem:trajectory-state} and its proof. The argument does not assume a bound on $1/(\alpha\gamma)$ or on the clock.

\textit{Step 1 (a forecast floor before attainment).}
On a pre-hit prefix put $r=\sqrt D$ and $V=(R+r)^2$. The radius estimates give $V<80R^2$, $3r(R+r)<192R^2$ and $\gamma>\varepsilon/[3r(R+r)]\ge\varepsilon/(3V)$. We claim
\begin{equation}\label{simp:floor}
\mathcal{Z}\ge\frac{\varepsilon^2}{288V\gamma},\qquad
\delta\ge\min\left\{\frac{\varepsilon^2}{2304V\gamma},\frac{D\gamma}{64(1+\chi)^2}\right\}.
\end{equation}
If $\sigma>\gamma/2$, the passing test bounds the certificate without its momentum term by $48\mathcal{Z}$. Completing the square in the affine model and testing $u=(1-t)y+tx^*$ gives
\[
\begin{aligned}
48\mathcal{Z}
&\ge F(y)-\inf_u\left\{F(u)+\frac{\sigma}{2}\norm{u-y}^2\right\}\\
&\ge \max_{0\le t\le1}\left\{t\varepsilon-\frac{\sigma V}{2}t^2\right\}
\ge\min\left\{\frac{\varepsilon}{2},\frac{\varepsilon^2}{2\sigma V}\right\}.
\end{aligned}
\]
The floor follows from $\sigma\le\gamma$ and $\varepsilon<3V\gamma$.
If $\sigma\le\gamma/2$, the stage began at a correction, since initially $\sigma=\gamma$ and a continuation from equality preserves it. On an uncorrected step with $M\le\sigma_{i-1}$, $\theta=1$ makes both weights equal before the common raise, so $\gamma=\sigma$ afterwards. Hence the part of this stage following its reset is clip-free. At the reset
\[
\frac{\mathcal{Z}_j}{\gamma_j-\sigma_j}\ge\frac{\norm d^2}{(\sigma_j)^2}+\frac{2e}{\sigma_j},\qquad\varepsilon<F(y_j)-F^*\le\sqrt V\norm d+e ,
\]
For $u=\norm d$, the smallest feasible error is bounded below by $[\varepsilon-\sqrt V u]_+$. Minimizing the resulting piecewise quadratic in $u\ge0$ gives
\[
 \frac{\mathcal Z_j}{\gamma_j-\sigma_j}\ge
 \begin{cases}
 \varepsilon^2/(V\sigma_j^2),&\varepsilon\le V\sigma_j,\\
 2\varepsilon/\sigma_j-V,&\varepsilon>V\sigma_j.
 \end{cases}
\]
The minimizer is $u=\min\{\varepsilon/\sqrt V,\sigma_j\sqrt V\}$. The bound decreases with $\sigma_j$.

Continuation has $\mathcal{Z}^+\ge(1-\theta)\mathcal{Z}$ and $(\theta^2/\alpha)-\sigma=(1-\theta)(\gamma-\sigma)$, and a raise multiplies numerator and denominator by the same factor and increases $\sigma$. So the bound holds at every later state of the stage with $\sigma$ in place of $\sigma_j$, also as $V$ grows. For $\sigma\le\gamma/2$ it gives $\mathcal{Z}\ge2\varepsilon^2/(V\gamma)$ if $\varepsilon\le V\sigma$, and $\mathcal{Z}\ge \varepsilon>\varepsilon^2/(3V\gamma)$ otherwise. Both imply~\eqref{simp:floor}.

\textit{Step 2 (the common floors).} The bound $\gamma>\varepsilon/(3V)$ and $V<80R^2$ give $\gamma>\varepsilon/(240R^2)$. If $\mu_f>0$, safety gives $\gamma\ge \sigma\ge\mu_f/2$. Hence
\begin{equation}\label{simp:gammafloor}
\gamma\ge\mu_{f,\varepsilon}/2
\end{equation}
on every pre-hit input. For $\mu_f>0$, combine the certificate lower bound from \cref{lem:energy-cert} with the passing test. Since $F(y)-F^*>\varepsilon$, this yields
\[
 \mathcal Z\ge\frac{\mu_f\sigma\varepsilon}{16(\mu_f+\sigma)(\sigma+\gamma)}
 \ge\frac{\mu_f\varepsilon}{96\gamma}.
\]
The last inequality uses $\sigma\le\gamma$ and $\mu_f\le2\sigma$. Together with~\eqref{simp:floor} and $V<80R^2$,
\begin{equation}\label{simp:commonfloor}
\mathcal{Z}\ge\max\Bigl\{\frac{\mu_f\varepsilon}{96\gamma},\ \frac{\varepsilon^2}{23040R^2\gamma}\Bigr\}=\frac{\mu_{f,\varepsilon}\varepsilon}{96\gamma},
\end{equation}
and for $\mu_f=0$ only the second term remains.

For the third inequality of \cref{lem:global-state}, the bound $F(y)-F^*\le3\gamma r(R+r)$ of \cref{lem:trajectory-state} and a gap above $\varepsilon$ give $\gamma r>\varepsilon/(3(R+r))$, so
$D\gamma^2>\varepsilon^2/(9V)>\varepsilon^2/(720R^2)$. Also $D\gamma^2\ge\mu_f^2D_0/4$ by $D\ge D_0$ and~\eqref{simp:gammafloor}. The larger of two numbers is at least their mean, so
\[
D\gamma^2\ge\frac18\Bigl(\mu_f^2D_0+\frac{\varepsilon^2}{180R^2}\Bigr) .
\]

\subsection{Proof of Lemma~\ref{lem:prefix-energy}}\label{app:shared}

For \cref{alg:main} assume $\mu_f>0$ and consider a finite prefix on which $1/\alpha\le M$ and $\mu_f/2\le \sigma\le\gamma\le M/2$, where $M\ge\mu_f$ is any bound for that prefix. Put $K=2M/\mu_f$, $A_K=2048K^3$, and $\sigma_0=1/(2\alpha_0)$. 

\textit{Step 1 (certificate and energy).} Put $K=2M/\mu_f\ge2$ and $\mathcal{P}=F(y)-F^*+\gamma\norm{v-y}^2/2$. The first inequality of \cref{lem:energy-cert} gives $\mathcal{P}\le(1+\sigma/\mu_f)\Cert_\sigma$. Every committed state passes the certificate, and $(1+\sigma/\mu_f)(1+\gamma/\sigma)=1+\sigma/\mu_f+\gamma/\sigma+\gamma/\mu_f\le2K$ by $\sigma\le\gamma\le M/2$ and $\sigma\ge\mu_f/2$, so
\begin{equation}\label{simple-comparison}
\mathcal{P}\le32K\mathcal{Z} .
\end{equation}

\textit{Step 2 (the cost of a correction).} Let $\mathcal P_{\rm old}$ be the input energy. The query bound and descent inequality give $\norm{x-y}^2/\alpha\le\gamma\norm{v-y}^2$ and $\norm{p-y}^2/(2\alpha)\le\mathcal P_{\rm old}+e$. By the squared triangle inequality and $v^+-p=(1-\theta)(p-y)/\theta$,
\[
 \alpha\norm d^2\le8\mathcal P_{\rm old}+4e,
 \qquad
 \frac{\theta^2}{\alpha}\norm{v^+-p}^2\le2(\mathcal P_{\rm old}+e).
\]
The raise has $r\le2$. At the reset, $\sigma\ge\mu_f/2$ and $\gamma\le M/2$. Substitution in the reset certificate gives
\[
 e+\frac{\alpha}{2}\norm d^2+\frac{\norm d^2}{2\sigma}
   +\frac\gamma2\norm{v^+-p}^2
 \le(6+4K)(\mathcal P_{\rm old}+e)\le7K(\mathcal P_{\rm old}+e).
\]
The reset forecast is at most $16K^2(\mathcal P_{\rm old}+e)$. By~\eqref{simple-comparison} and $e\le\mathcal Z_{\rm old}/8$, it is at most $(512K^3+2K^2)\mathcal Z_{\rm old}\le520K^3\mathcal Z_{\rm old}$.

\textit{Step 3 (forecast).} Before raising, the allowance $e\le\theta\mathcal Z/8$ bounds the forecast by $(1-3\theta/4)\mathcal Z$. A correction costs at most $520K^3$, and $520e^{3/4}<2048$, so multiplying the step estimates gives
\begin{equation}\label{simple-global-ledger}
\mathcal{Z}_n\le \mathcal{Z}_0\exp(-3\chi_n/4+U_n)A_K^{J_n}.
\end{equation}
\looseness=-1 Here $U_n$ includes raises on correction steps, and $J_n$ counts all corrections. Combining~\eqref{simple-comparison} with $J_n\le[\log(\sigma_0/\sigma_n)+U_n-V_n]/\log2$ from~\eqref{mastercount} and $\sigma_n\ge\mu_f/2$ gives
\begin{equation}\label{masterrate}
  \mathcal{P}_n\le32K\mathcal{Z}_0\exp(-3\chi_n/4+U_n)\,A_K^{[\log(\sigma_0/(\mu_f/2))+U_n-V_n]/\log2}.
\end{equation}

\subsection{Proof of Lemma~\ref{lem:global-clock-main}}\label{app:clock-budget}

Fix a pre-hit prefix satisfying $1/\alpha_i\le Q\gamma_{i-1}$, put $U=U_n$ and $\mathcal L=\log(eQR^2/(\alpha_0\varepsilon))$, and let $C$ denote absolute constants. As $\mu_{f,\varepsilon}\ge\varepsilon/(240R^2)$ and $R^2/(\alpha_0\varepsilon)>2$, both $\log Q$ and $1+\log(1/(\alpha_0\mu_{f,\varepsilon}))$ are $O(\mathcal L)$.

\textit{Step 1 (the two clock regimes).} On inputs with $\sigma\le\gamma/2$, the normalized weight $e^{-U_i}\gamma_i$ contracts by a factor at most $1-\theta_i/2$. It starts at $1/(2\alpha_0)$, is nonincreasing on the other steps, and stays above $e^{-U}\mu_{f,\varepsilon}/2$. Thus these inputs spend at most $C(\mathcal L+U)$ clock.

It remains to count the clock spent with $\sigma>\gamma/2$. We first count the stages that reach this regime and then bound the clock spent in each. A stage starts at initialization or a correction. On continuation $e^{-U_i}\sigma_i$ is nonincreasing, and a correction reduces it by at least two. Any stage reaching $\sigma>\gamma/2$ must have normalized modulus above $e^{-U}\mu_{f,\varepsilon}/4$. There are therefore at most $C(\mathcal L+U)$ such stages. Within a stage $(\gamma-\sigma)/\sigma$ is nonincreasing, dropping to zero at a strict clip, so it enters this regime at most once.

To bound the clock within one such stage, write $\Gamma=e^U/(2\alpha_0)$ for an upper bound on every weight. At the start $j$ of a stage reaching $\sigma>\gamma/2$ we have $\sigma_j>e^{-U}\mu_{f,\varepsilon}/4$. At a reset, \cref{lem:trajectory-state}, the allowance~\eqref{allow}, and the proximal inequalities in \cref{app:step} give $\mathcal P_{\rm in}+e\le C\Gamma R^2$, $\|d\|^2\le CQ\Gamma^2 R^2$, and $(\theta^2/\alpha)\|v^+-p\|^2\le C\Gamma R^2$. Substitution in the reset formula gives
\[
 \mathcal Z_j\le CQ\Gamma R^2(\Gamma/\sigma_j)^2,\qquad
 e^{-U_i}\mathcal Z_i\ge e^{-U}\mu_{f,\varepsilon}\varepsilon/(96\Gamma)
\]
where the second inequality is~\eqref{simp:commonfloor}. The first also holds at initialization since $D_0\le4R^2$. Dividing these bounds shows that the logarithm of the ratio is at most
$\log(CQR^2/(\alpha_0\varepsilon))+3\log(1/(\alpha_0\mu_{f,\varepsilon}))+7U=O(\mathcal L+U)$.
On continuation $e^{-U_i}\mathcal Z_i$ contracts by at least the fraction $1-e^{-3\theta_i/4}$. Each such stage therefore spends $O(\mathcal L+U)$ clock, including its final correction. Multiplication by the stage count gives

\begin{equation}\label{simp:clock}
 \chi\le C\bigl(\mathcal L+U\bigr)^2.
\end{equation}

\textit{Step 2 (the strongly convex clock inequality).}
Suppose $\mu_f>0$. On any pre-hit prefix satisfying the ratio bound, normalized-weight monotonicity gives $\gamma\le e^U/(2\alpha_0)$. \Cref{lem:prefix-energy} and the reset estimates in Appendix~\ref{app:shared} apply with $K=2Qe^U/(\alpha_0\mu_f)$. Use $\mathcal{Z}_0\le24R^2/\alpha_0$, $J_n\le[\log(1/(\alpha_0\mu_f))+U_n]/\log2$, and $\varepsilon<\mathcal{P}_n\le32K\mathcal{Z}_n$ in~\eqref{simple-global-ledger} of Appendix~\ref{app:shared}. Taking logarithms yields
\[
 \chi_n\le C\left[\log\frac{eR^2}{\alpha_0\varepsilon}
             +\log^2\frac{eQ\alpha_0^{-1}}{\mu_f}+U_n^2\right].
\]

\textit{Step 3 (closing the bound).}
Put $A_0=\log(eR^2/(\alpha_0\varepsilon))+\log^2(eQ/(\alpha_0\mu_{f,\varepsilon}))$. If $\mu_{f,\varepsilon}=\mu_f$, Step 2 gives $\chi_n\le C(A_0+U_n^2)$. Otherwise $\mu_{f,\varepsilon}=\varepsilon/(240R^2)$ and $\mathcal L^2\le A_0$, so~\eqref{simp:clock} gives the same inequality. Since $U_n\le\log(1+\chi_n)$, apply $C\log^2 x\le x/2+B_C$ with $x=1+\chi_n$. Rearranging gives
\[
 1+\chi_n\le 2(1+B_C+CA_0)<C_0A_0/2 .
\]
Here $A_0\ge1$, so a single absolute $C_0$ suffices.

\subsection{Proof of Theorem~\ref{bi:horizon}}\label{bi:proof-horizon}
The initial ratio is two and \cref{lem:regional-main} supplies the initial floor. Assume the floor and the clock bound hold on the previous pre-hit prefix. For a regional upper bound, proximal nonexpansiveness on the acceptance interval gives
\[
 \|p-x^*\|\le5R+\alpha G\le10R.
\]
For a global upper bound, the upper model applies directly. In either case,~\eqref{bi:upper}, $1/(2\alpha)\ge Qq\lambda/8$ and the tolerance from \cref{lem:regional-main} and the root bound give the upper test. On the same interval $M\ge Q\gamma/2\ge2\gamma$ also pays for the full prospective raise. The inheritance and acceptance arguments in \cref{app:ratio-progress,app:allowance-test} keep the ratio within $Q$, including at the hitting step.

If the enlarged prefix is still pre-hit, \cref{app:clock-budget} with $\mu_f=0$ gives $1+\chi_k<\tau$ before the floor at $k$ is used. \Cref{lem:regional-main} then gives the floor at $k$, closing the induction. Since each increment exceeds $1/(2\sqrt Q)$, the clock bounds the number of accepted steps. The trial count in \cref{app:ratio-progress} then gives the oracle bound.

An admissible ratio exists: after subtracting the tangent slope at $x^*$, bounded gradients on $T$ give $\widehat\mu_T(t)\le Mt$ for a finite $M$. Completing the square gives~\eqref{bi:upper} when $Q\lambda s\ge2M^2$, eventually true since $Qs\to\infty$ and $\lambda\ge m_S(\varepsilon/(16\cdot64^8))/16>0$.

\subsection{Proof of Corollary~\ref{bi:power}}\label{bi:proof-power}
We first bound the ratio needed for acceptance, check localization, and then insert the regional lower modulus. Consider a power contribution $c_\nu L_\nu t^{1+\nu}$ to the upper model. Set $u=qt$ and multiply \eqref{bi:upper} by $q$. The power term acquires a factor $q^{-\nu}\le1$. Even if only half the quadratic budget is available, it suffices to check $c_\nu L_\nu u_0^{1+\nu}\le s$ at $u_0=\sqrt{16s/(Q\lambda)}$. Below $u_0$ the allowance suffices. Above $u_0$ the quadratic term grows at least as fast as the power term. Substituting the regional allowance gives the sufficient condition
\begin{equation}\label{eq:regional-power-test}
 Q^{(1+3\nu)/4}\ge C_\nu L_\nu\lambda^{-(1+\nu)/2}\varepsilon^{-(1-\nu)/2}\tau^{1-\nu}.
\end{equation}

For this corollary use $\lambda=m_S(\varepsilon/(16\cdot64^8))/16$, which is no larger than the scale in the theorem. The H\"older tangent bounds give $c_\nu=2^{1-\nu}/(1+\nu)$ and $G\le L_\nu(5R)^\nu$. Put $\zeta=L_\nu R^{1+\nu}/\varepsilon$. The diameter bound and the chord of length $R$ in $S$ give $c\varepsilon/R^2\le\lambda\le C_\nu(1+\zeta)\varepsilon/R^2$. Hence $L_\nu R^{\nu-1}/\lambda=O_\nu(1)$ if $\zeta\le1$. For $\zeta\ge1$, the upper bound on $\lambda$ gives $Q_\nu\ge c_\nu' L_\nu R^{\nu-1}/\lambda$, since
\[Q_\nu/(L_\nu R^{\nu-1}/\lambda)=(\zeta^3\varepsilon/(\lambda R^2))^{(1-\nu)/(1+3\nu)},\]
where $Q_\nu=(L_\nu^2/(\lambda^{1+\nu}\varepsilon^{1-\nu}))^{2/(1+3\nu)}$. 

For both this proof and the mixed-model proof below, let $Q_0\ge\max\{(1-\rho)^{-2},Q_\nu\}$, set $L_Q=\log(eQ_0(1+R^2/(\alpha_0\varepsilon)))$ and take $Q=D_\nu Q_0 L_Q^{8(1-\nu)/(1+3\nu)}$. Then $\tau\le C_\nu(1+\log^2 D_\nu)L_Q^2$. The powers of $L_Q$ cancel in \eqref{eq:regional-power-test}, and a sufficiently large fixed $D_\nu$ absorbs the remaining logarithm. Here $Q_0=\max\{(1-\rho)^{-2},Q_\nu\}$ also covers localization. Inserting both lower bounds for $\lambda$ from \cref{bi:model} gives the displayed exponents. In the power-law branch the replacement of $\varepsilon/16$ by $\varepsilon/(16\cdot64^8)$ contributes a fixed factor to the bound on $1/\lambda$. The exponent of accuracy is unchanged.

For \cref{prop:mixed} with $\nu=1$, $Q=\max\{(1-\rho)^{-2},4(L_2+L_1)/\mu_f\}$ suffices in \cref{thm:master}. Otherwise averaging~\eqref{eq:mixed-bound} at the barycenter gives $\widehat\mu(t)\le L_2t^2/2+c_\nu L_\nu t^{1+\nu}$. Take $\lambda_*=\mu_f/16\le m_S(s_0)/16$. Half the quadratic budget covers $L_2t^2/2$ when $Q\lambda_*\ge8L_2$. Apply the preceding power calculation at $\lambda_*$ with $Q_0=\max\{(1-\rho)^{-2},L_2/\mu_f,Q_\nu\}$ and enlarge $D_\nu$ to cover both terms. The logarithms have the stated form since $R^2\le2(F(z_0)-F^*)/\mu_f$ and, on arbitrarily long chords,~\eqref{eq:mixed-bound} forces $L_2\ge\mu_f$.

\subsection{Proof of Theorem~\ref{thm:master}}\label{app:global-horizon}
Fix $0<\varepsilon<F(p_0)-F^*$ and use the constant from \cref{app:clock-budget}.

\textit{Step 1 (existence of an admissible ratio).}
For affine $f$, condition~\eqref{masteruniv} is automatic. Otherwise fix $t_1>0$ and put $s_1=\widehat{\mu}(t_1)>0$. As $Q\to\infty$, the allowance in~\eqref{eq:horizon-allowance} has order $Q^{-1/2}(\log Q)^{-4}$, so $s\to0$ and $Qs\to\infty$. Eventually $s\le s_1$ and
\[
 \sqrt{\frac{8s}{Q\mu_{f,\varepsilon}}}\le\frac{s}{s_1}t_1.
\]
Monotonicity and chord scaling give~\eqref{masteruniv}, proving part (i).

\textit{Step 2 (propagating the ratio and clock bounds).}
Fix an admissible $Q$ and put $\tau=\tau(Q)$. At initialization the ratio is two and $1+\chi_0<\tau/2$. Suppose the previous ratios are at most $Q$ and the current point is pre-hit with $1+\chi<\tau/2$. The state bounds in \cref{lem:global-state} control both branches of the allowance.

Put $q=2\gamma/\mu_{f,\varepsilon}\ge1$. The forecast bound gives $\mathcal{Z}/8\ge\varepsilon/(384q)$. The radius bound and $1+\chi\le\tau$ give
\[
 \frac{D\gamma}{64(1+\chi)^2}
 \ge\frac{\mu_f^2D_0+\varepsilon^2/(180R^2)}{512\gamma\tau^2}
 =\frac1q\frac{\mu_f^2D_0+\varepsilon^2/(180R^2)}{256\mu_{f,\varepsilon}\tau^2}.
\]
Consequently,
\[
 \delta\ge\frac1{384q}
 \min\left\{\varepsilon,
 \frac{\mu_f^2D_0+\varepsilon^2/(180R^2)}{\mu_{f,\varepsilon}\tau^2}\right\}.
\]
By \cref{app:ratio-progress}, $\theta\ge(\sqrt{17}-1)/(4\sqrt Q)$. Since $(\sqrt{17}-1)/1536>1/1300$, multiplication of the preceding inequality gives $\theta\delta\ge s/q$ for $s$ equal to the displayed minimum divided by $1300\sqrt Q$.

Inheritance (\cref{lem:progress-main}) and the upper test (Appendix~\ref{app:allowance-test}) then give an accepted step with
\[
 1/\alpha\le Q\gamma_{\rm in},\qquad \theta>1/(2\sqrt Q).
\]
If the new point is still pre-hit, \cref{lem:global-clock-main} gives $1+\chi<\tau/2$ again. This closes the induction and bounds the ratio through the hitting step.

\textit{Step 3 (iteration and oracle counts).}
If $N=\lceil2\sqrt Q\,\tau\rceil$ steps were all pre-hit, their clock would exceed $N/(2\sqrt Q)\ge\tau$, contrary to Step~2. This proves~\eqref{eq:horizon-count}. Since $\alpha_0/\alpha_n\le Qe^{U_{n-1}}/2$ and $1+\chi_{n-1}\le\tau$, \cref{app:ratio-progress} bounds rejections by $n+\log_2(Q\tau)$. Adding initialization and accepted trials gives the oracle bound. A stationary return already solves the problem.

\subsection{\texorpdfstring{Proof of \cref{cor:convex}}{Proof of the convex corollary}}\label{app:convex}
For $\mu_f=0$, $\mu_{f,\varepsilon}=\varepsilon/(240R^2)$ and the second numerator in~\eqref{eq:horizon-allowance}, divided by $\mu_{f,\varepsilon}$, equals $4\varepsilon/3$. We verify~\eqref{masteruniv} first for power bounds and then for the intrinsic gauge. Oracle counts add only a constant factor and a logarithmic term.

\textit{Step 1 (H\"older gradients).}
Put $q=\max\{(1-\rho)^{-2},(L_\nu R^{1+\nu}/\varepsilon)^{4/(1+3\nu)}\}$, $L_Q=\log(eqR^2/(\alpha_0\varepsilon))$ and $\eta=(1-\nu)/(1+3\nu)$. At $Q=D_\nu qL_Q^{8\eta}$ the clock is $O((1+\log^2D_\nu)L_Q^2)$ and $s\gtrsim\varepsilon/(\sqrt Q\,\tau^2)$. The power bound thus gives
\[
 \frac{\widehat\mu(\sqrt{8s/(Q\mu_{f,\varepsilon})})}{s}
 \le C_\nu\frac{L_\nu R^{1+\nu}}{\varepsilon}Q^{-(1+3\nu)/4}\tau^{1-\nu}
 \le C'_\nu D_\nu^{-(1+3\nu)/4}(1+\log^2D_\nu)^{1-\nu}.
\]
Choose $D_\nu$ so the last term is at most one. The cost is $O_\nu(\sqrt q\,L_Q^{2+4\eta})$, including the endpoints. Since $R^2/(\alpha_0\varepsilon)>2$, $L_Q=O_{\nu,\rho}(\log((R^2/\alpha_0+L_\nu R^{1+\nu})/\varepsilon))$.

\textit{Step 2 (an arbitrary finite curvature bound).}
The chord inequality gives $k\widehat{\mu}(2R/k^{3/2})\le\widehat{\mu}(2R)/\sqrt k\to0$, so the gauge $k=\mathcal N_\varepsilon$ is finite. Let
\[
 q=\max\{(1-\rho)^{-2},k^2\},\qquad
 L_Q=\log\frac{eq\alpha_0^{-1}R^2}{\varepsilon},\qquad
 Q=Dq L_Q^8,
\]
with a sufficiently large absolute $D$. As in Step 1,
$s\ge c\varepsilon/[\sqrt Q(1+\log^2D)^2 L_Q^4]$, while the gauge gives $\widehat{\mu}(2R/k^{3/2})\le\varepsilon/k$.
Set $t=\sqrt{8s/(Q\mu_{f,\varepsilon})}$. Then
\[
 \left(\frac{t}{2R/k^{3/2}}\right)^2
       =\frac{480sk^3}{Q\varepsilon},\qquad
 Qs\ge\frac{c\sqrt D\sqrt q\,\varepsilon}{(1+\log^2D)^2}.
\]
Choosing $D$ so that $Qs\ge480k\varepsilon$, and using
$s\le\varepsilon/(1300\sqrt Q)$, also ensures $t\le2R/k^{3/2}$. Chord scaling now gives
$\widehat{\mu}(t)\le\frac{t}{2R/k^{3/2}}\frac{\varepsilon}{k}\le s.$
This proves feasibility. The horizon is $O(\sqrt q\,L_Q^6)$, which is
$O_\rho(k\log^6(\alpha_0^{-1}R^2k/\varepsilon))$. Adding the initial search proves part (ii).

\subsection{\texorpdfstring{Proof of \cref{prop:smooth-energy}}{Proof of the smooth energy guarantee}}\label{app:strong}

Assume that $f$ is $\mu_f$-strongly convex with an $L$-Lipschitz gradient. Put $M=\max\{1/\alpha_0,4L\}$ and $K=2M/\mu_f$. By \cref{prop:basic}, $\sigma\ge\mu_f/2$.

\textit{Step 1 (scale bounds).} Initially $1/\alpha\le M$ and $\sigma\le\gamma\le M/2$. Both trial tests pass whenever $1/\alpha\ge L$, because $c\le L$ and $e=0$. Hence $1/\alpha_k\le\max\{\rho/\alpha_{k-1},2L\}\le M$, and the probe initialization gives $1/\alpha_k\le2L$ by induction. Without a raise, $\gamma^+=\theta^2/\alpha\le\gamma$. Otherwise $r\le1+\theta$ and $r\sigma_{\rm tr}\le c\le L\le M/4$ give
\[
 \gamma^+\le
 \begin{cases}
 (1+\theta)M\theta^2\le3M/8,&\theta\le1/2,\\
 (1-\theta^2)M/2+\theta M/4\le M/2,&\theta>1/2.
 \end{cases}
\]
The second branch uses the root equation. The raise preserves $\sigma\le\gamma$, and a correction only decreases $\sigma$.

\textit{Step 2 (energy at every iteration).} The scale bounds allow Steps 1--3 of \cref{app:shared} to be applied to every prefix, so~\eqref{masterrate} holds.
Since $\sigma_{\rm tr}\ge\mu_f/2$ (on a strict clip $M\ge c\ge\mu_f$), $\theta^2/\alpha\ge\mu_f/2$ and every accepted step adds at least $1/\sqrt K$ to the clock. For the probe initialization $K=8\kappa$, $1/(\alpha_0\mu_f)<2\kappa$ and $\mathcal Z_0\le12(F(z_0)-F^*)$. Taking logarithms in~\eqref{masterrate} and using $U_n\le\log(1+\chi_n)$ and $V_n\ge0$ gives
\[
 \log\frac{\mathcal P_n}{F(z_0)-F^*}
 \le -3\chi_n/4+C\log^2(2\kappa)+C\log(2\kappa)\log(1+\chi_n).
\]
For $b\ge1$, maximizing $b\log(1+t)-t/8$ over $t\ge0$ bounds it by $b\log(8b)$. Thus the last term costs at most $\chi_n/8+C'\log^2(2\kappa)$. With the probe, $M_k\le2L$ gives $\chi_n\ge n/(2\sqrt\kappa)$. Before absorbing the logarithmic term, division by $n$ gives the stated asymptotic exponent. Using the absorbed bound instead proves~\eqref{polyrate}, since the accuracy logarithm absorbs an absolute constant even at $\kappa=1$.

For $n$ accepted steps, the rejection count telescopes to $n\log_2(1/\rho)+\log_2(\alpha_0/\alpha_n)\le n+\log_2(2\kappa)$ by the scale bounds. Each trial uses at most two calls per oracle, giving the first-attainment count.

\textit{Step 3 (convergence).} By~\eqref{polyrate} $\mathcal{P}_n\to0$, so $y_n\to x^*$ by strong convexity, and $\gamma_n\ge\mu_f/2$ gives $v_n-y_n\to0$. The query points are convex combinations of these points and converge to $x^*$ as well.

\textit{Step 4 (trajectory ball).} By~\cref{lem:trajectory-state}, $z_0$, all iteration queries, primal and momentum points, and cached anchors lie in $\mathcal K$. The accepted secants and $1/\alpha_0$ are therefore at least $\mu_{\mathcal K}$. The lower models used in \cref{app:basic,lem:energy-cert} involve only these points and $x^*$, so pre-attainment corrections still start above $\mu_{\mathcal K}$ and return at least $\mu_{\mathcal K}/2$. For~\eqref{eq:initial-choice}, $z_0$ is interior to $\mathcal K$ when $R>0$. The initial part of a nonzero probe therefore has strictly positive directional curvature, and convexity along the remaining segment gives a positive probe secant. Together with global smoothness this yields \[\mu_{\mathcal K}\le1/\alpha_0<2L.\] The probe endpoint may lie outside $\mathcal K$, so its secant need only be positive. Here $N_{\rm init}$ retains the actual initialization cost. Steps 1--3 and~\eqref{mastercount} apply with $\mu_{\mathcal K}$ in place of $\mu_f$.

\section{\texorpdfstring{Proof of \cref{thm:simple-local-recovery}}{Proof of the recovery theorem}}\label{app:simple-local-recovery}

\textit{Step 1 (localized safety).} Start at $j$ from part (i). For $k>j$, every tested endpoint and the query belong to the local neighborhood. Each secant is at least $\mu_\ell$, hence $b_k\ge\mu_\ell$. A strict clip gives $M\ge c\ge\mu_\ell$ before correction. Equality $M=\mu_\ell$ on a strict clip gives $F(p)=F^*$ by the completion of the square in Appendix~\ref{app:basic} and $v^+=p$, so the next query returns a minimizer. This is excluded on an infinite run. After a local reset, the certificate induction of Appendix~\ref{app:basic} uses only points in the neighborhood. Thus only the stage inherited at $j$ can correct a tentative estimate below $\mu_\ell$.

\textit{Step 2 (recovery and its cost).} The raise cap and corrections make $\sigma_k/(1+\chi_k)$ nonincreasing, proving part (ii). By Step~1, the floor $\mu_\ell/2$ becomes permanent at $\mu_\ell$ in the inherited stage or at $\mu_\ell/2$ in a local stage.

A strict clip gives $M\ge\mu_\ell$ before testing and its correction, if any, starts a local stage at or above the permanent floor. Until that floor is secured, at most one correction halves the estimate. Every raise takes the full clock factor, since a secant-limited raise would reach $c_k\ge\mu_\ell$. Telescoping yields
\[
 \sigma_n\ge\frac{\sigma_{j}}{2}\frac{1+\chi_n}{1+\chi_{j}}
 \qquad\text{before recovery}.
\]
Recovery occurs by the time this bound reaches $\mu_\ell$. Throughout the tail $\sigma_k\ge m_j$. Indeed, the inherited stage can cut an estimate below $\mu_\ell$ only once, and each later correction returns at least $\mu_\ell/2$. A strict clip has $M\ge c\ge\mu_\ell$, so also $\sigma_{\rm tr}\ge m_j$. The scale recurrence in Step~1 of Appendix~\ref{app:strong} gives $M_k\le H_j$ for $k>j$. Thus $\theta_k^2=w_k/M_k\ge m_j/H_j$ at every step after $j$. This turns the clock bound into (i), with the sharper increment $1/(2\sqrt\kappa)$ under the hypotheses of \cref{prop:smooth-energy}. From $k_0$ onward $M_n\le\max\{M_{k_0},2L\}$, so the trial count of Appendix~\ref{app:ratio-progress} gives the asserted oracle cost.

\textit{Step 3 (energy after recovery).} Fix $N$ after recovery. The scale induction of Appendix~\ref{app:strong}, Step 1, uses only global smoothness to bound $M_i$ and $\gamma_i$. On this tail $\sigma_i\ge\mu_\ell/2$. The certificate-to-energy estimate in Appendix~\ref{app:energy-cert} uses strong convexity only on $[y_i,x^*]$, which lies in the local neighborhood. Thus the per-step forecast and reset estimates of Appendix~\ref{app:shared} apply from state $N$ with $\mu_\ell$ and a fixed scale bound. By~\eqref{mastercount}, $U_n,J_n=O(\log(1+n))$ since $\sigma_n$ has a positive floor. Multiplication from $N$ gives
\[
 \log\mathcal P_n\le -3(\chi_n-\chi_N)/4+O(\log(1+n)).
\]
The clock lower bound of Step 2 is linear in $n-N$. Thus $\mathcal P_n\to0$, local strong convexity gives $y_n\to x^*$, and the weight floor gives $v_n\to x^*$ as well.

\textit{Step 4 (late scales).} Every trial query is a convex combination of the incoming primal and momentum points, so all trial queries converge uniformly to $x^*$. Proximal nonexpansiveness bounds every trial by $\|p_\alpha-x^*\|\le\|x-x^*\|+\alpha\|\nabla f(x)-\nabla f(x^*)\|$. Since every trial after recovery has $M\ge M_{k-1}/2\ge\sigma_{k-1}/2\ge\mu_\ell/4$, every sufficiently late trial proximal point is local as well. Thus $M_k\le\max\{\rho M_{k-1},2L_\ell\}$, giving $\limsup_k M_k\le2L_\ell$. The recovered floor then yields $\liminf_k\theta_k\ge\sqrt{\mu_\ell/(4L_\ell)}$. If a late estimate is at most $L_\ell$, the secant cap keeps subsequent estimates at most $L_\ell$. For $k>r$ on this tail, the raise cap and root equation give
\[
 \gamma_k\le (1-a_r\theta_k)\gamma_{k-1}+\theta_k L_\ell,
 \qquad a_r=\chi_r/(1+\chi_r).
\]
Iteration gives $\limsup_k\gamma_k\le L_\ell/a_r$ for each sufficiently late $r$, since $\sum_k\theta_k=\infty$. Here $a_r$ is positive for $r\ge1$. Sending $r\to\infty$ gives the conditional weight bound. Independently of that condition, the lower bound on $\liminf_k\theta_k$ and the tail energy inequality of Step 3 give the exponent in (iii).

\end{document}